\documentclass{amsart}

\usepackage{amsfonts}
\usepackage{amssymb}
\usepackage{latexsym}
\usepackage{diagrams}
\newarrow{Mapsto}{|}{-}{-}{-}{>}

\newtheorem{theorem}{Theorem}[section]
\newtheorem{lemma}[theorem]{Lemma}
\newtheorem{definition}[theorem]{Definition}
\newtheorem{proposition}[theorem]{Proposition}

\newtheorem{corollary}[theorem]{Corollary}

\theoremstyle{definition}

\newtheorem{example}[theorem]{Example}

\numberwithin{equation}{section}

\newcommand{\R}{\mathbb{R}}
\newcommand{\Z}{\mathbb{Z}}
\newcommand{\F}{\mathbb{F}}
\newcommand{\Imt}{\mbox{Im}\,}

\begin{document}

\title
{Generalized twisted sectors of orbifolds}

\author{Carla Farsi}
\address{Department of Mathematics, University of Colorado at Boulder, Campus
Box 395, Boulder, CO 80309-0395 } \email{farsi@euclid.colorado.edu}

\author{Christopher Seaton}
\thanks{The second author was partially supported by a Rhodes College Faculty Development Endowment Grant.}
\address{Department of Mathematics and Computer Science,
Rhodes College, 2000 N. Parkway, Memphis, TN 38112}
\email{seatonc@rhodes.edu}



\begin{abstract}

For a finitely generated discrete group $\Gamma$, the $\Gamma$-sectors of
an orbifold $Q$ are a disjoint union of orbifolds
corresponding to homomorphisms from $\Gamma$ into a groupoid presenting $Q$.
Here, we show that the inertia orbifold and $k$-multi-sectors are special cases
of the $\Gamma$-sectors, and that the $\Gamma$-sectors are
orbifold covers of Leida's fixed-point sectors. In the case of a
global quotient, we show that the $\Gamma$-sectors correspond to
orbifolds considered by other authors for global quotient orbifolds as well as their
direct generalization to the case of an orbifold given by a quotient by a Lie group.
Furthermore, we develop a model for the $\Gamma$-sectors corresponding to a
generalized loop space.

\end{abstract}

\maketitle


\section{Introduction}
\label{sec-intro}

In \cite{farsiseaton1}, the authors introduced the $\Gamma$-sectors
of an orbifold in order to determine a complete obstruction to the
existence of a nonvanishing vector.  The definitions of these
sectors was heavily motivated by several existing constructions for
orbifolds by Kawasaki (\cite{kawasaki1}, \cite{kawasaki2}, and
\cite{kawasaki3}), Chen and Ruan (\cite{chenruanorbcohom},
\cite{ruansgt}), Bryan and Fulman
(\cite{bryanfulman}), and Tamanoi (\cite{tamanoi1} and
\cite{tamanoi2}).

The goal of this paper is to show explicitly how the
$\Gamma$-sectors generalize these constructions.  In particular, we
show that the inertia orbifold corresponds to the $\Z$-sectors and
the $k$-multi-sectors correspond to the $\F_k$-sectors where $\F_k$
is the free group with $k$ generators. The orbifolds whose Euler
characteristics are considered by Bryan-Fulman and Tamanoi for
global quotients correspond to the $\Z^k$-sectors and
$\Gamma$-sectors, respectively, in the case that $Q$ can be
expressed as a global quotient; i.e. a quotient of a manifold by a
finite group.  Additionally, we show that the fixed-point sectors
introduced by Leida in \cite{leida} are orbifold-covered by the
$\Gamma$-sectors for an appropriate choice of $\Gamma$.

The work of Lupercio and Uribe in \cite{luperciouribeloop} (see also
\cite{defernetluperciouribe}) demonstrates that the inertia orbifold
naturally appears when considering the loop space of an orbifold.
Here, we show that the same holds true for the $\Gamma$-sectors; in
particular, they appear when considering smooth maps $M_\Gamma
\rightarrow Q$ where $M_\Gamma$ is a smooth manifold with
fundamental group $\Gamma$.  This generalizes results of Tamanoi
in \cite{tamanoi2}, stated for global quotients in the context of
orbifold bundles.

In the case that an orbifold $Q$ is presented by a quotient $M/G$
where $M$ is a manifold and $G$ is a Lie group acting locally
freely, i.e. properly with discrete stabilizers, there is a very
natural extension of the definition of the orbifolds considered by
Bryan-Fulman and Tamanoi (see Definition
\ref{def-globalgammasectors}).  We show that this again coincides
with the $\Gamma$-sectors.  Note, however, that such a presentation
of the $\Gamma$-sectors leads to a different indexing of the
sectors.  In the case of a global quotient, the $\Gamma$-sectors are
naturally indexed by $G$-conjugacy classes $(\phi)_\sim$ of
homomorphisms $\phi : \Gamma \rightarrow G$ whose images fix a
nonempty subset of $M$; we use $t_{M;G}^\Gamma$ to denote the set of
conjugacy classes of such homomorphisms (see Subsection
\ref{subsec-globaldef}).  On the other hand, if $\mathcal{G}$ is an
orbifold groupoid presenting $Q$, then the sectors are indexed by
elements of $T_Q^\Gamma$, the set of $\approx$-classes of elements
of $\mbox{HOM}(\Gamma, \mathcal{G})$ or equivalently connected
components of $|\mathcal{G} \ltimes \mbox{HOM}(\Gamma,
\mathcal{G})|$ (see \cite[Subsection 2.2]{farsiseaton1} or
Subsection \ref{subsec-localdef} below).  The discrepancy arises
from the fact that the fixed-point set of a homomorphism $\phi :
\Gamma \rightarrow G$ need not be connected, and hence the
$\Gamma$-sector corresponding to $(\phi)_\sim \in t_{M;G}^\Gamma$
may correspond to the disjoint union of several $\Gamma$-sectors,
each corresponding to one element of $T_Q^\Gamma$ (see Example
\ref{ex-localglobaldiffer}).

Note that in \cite{farsiseaton1}, we made the requirement that our
local groups act with a fixed-point set of codimension 2; however,
it was noted that the construction of the $\Gamma$-sectors did not
require this property. In this paper, we do not make this
requirement.  Moreover, while our primary interest is the case
of orbifold groupoids, many of these constructions and results generalize directly
to the case of orbispaces (see \cite{chenorbispace}).
We mention these generalizations as they arise.

The first author would like to thank the MSRI for its hospitality
during the preparation of this manuscript.


\section{Two Definitions of the $\Gamma$-Sectors for Quotient Orbifolds}
\label{sec-definitions}

In \cite{farsiseaton1}, the $\Gamma$-sectors
of a general orbifold were constructed in terms of the orbifold structure given by an
orbifold groupoid $\mathcal{G}$; that is, a proper, \'{e}tale Lie
groupoid.  For background on orbifolds from this perspective,
the reader is referred to \cite{ademleidaruan};
see also \cite{moerdijkmrcun} and \cite{moerdijkorbgroupintro}.
In Subsection \ref{subsec-globaldef}, we are primarily
concerned with orbifolds presented as the quotient of
a manifold by a Lie group.  We construct the $\Gamma$-sectors
directly from such a presentation. This construction was introduced
by Tamanoi in \cite{tamanoi1} and \cite{tamanoi2} for the case that
$G$ is finite; the definitions are unchanged for general $G$.  In Subsection
\ref{subsec-localdef}, we review the key points of the construction
in \cite{farsiseaton1} and give other interpretations.  Note that we
use slightly different notation for the $\Gamma$-sectors of a
quotient orbifold to distinguish from the construction using a
general orbifold groupoid; these definitions will be compared in
Section \ref{sec-localvsglobaldef}.


\subsection{$\Gamma$-Sectors of a Quotient Presentation}
\label{subsec-globaldef}

Let $Q$ be an $n$-dimensional \emph{quotient orbifold}.  By this,
we mean that $Q$ is presented by $G \ltimes M$ where $M$ is a smooth
manifold, $G$ is a Lie group acting
smoothly on $M$, and $G\ltimes M$ is Morita equivalent to an
orbifold groupoid, i.e. a proper \'{e}tale Lie groupoid.  In
\cite[page 536]{ademruan} and \cite[page 57]{ademleidaruan} (see
also \cite[page 76]{kawasaki1}), it is noted that this is the case
whenever the following conditions are satisfied:
\newcounter{Fcount}
\begin{list}{\roman{Fcount}.}{\usecounter{Fcount}}
\item   the isotropy group $G_x$ for each $x \in M$ is finite,
\item   there is a smooth slice $S_x$ at each $x \in M$, and
\item   for each $x, y \in M$ with $y \notin Gx$, there are slices $S_x$
and $S_y$ such that $GS_x \cap GS_y = \emptyset$.
\end{list}
In particular, it is noted that (ii) and (iii) are automatically
satisfied if $G$ is compact.  The following special cases are worth
noting; occasionally, we will restrict our attention to one of
these.
\begin{itemize}
\item   If $G$ is a finite group, then $Q$ is a {\bf global quotient
        orbifold}.
\item   If $G$ is a discrete group acting properly discontinuously,
        then $Q$ is a {\bf good} orbifold (see
        \cite[Definition 13.2.3]{thurston} or \cite[page 20]{boileaump}).
\end{itemize}
Note that we use the notation $M/G$ to indicate the quotient space
as a topological space only; the orbifold (i.e. Morita equivalence class
of the groupoid $G \ltimes M$) will generally be denoted $Q$.
In \cite{henriquesmetzler}, the question of whether every orbifold
can be expressed as a quotient is addressed.  In general, this
question remains unresolved.

Note that in the case of a good orbifold (including the case of a
global quotient), the groupoid $G \ltimes M$ is an orbifold
groupoid. On the other hand, if $G$ is a Lie group of positive
dimension, then $G \ltimes M$ is not \'{e}tale, though it is Morita equivalent to an orbifold
groupoid.  In general, $G \ltimes M$ as well
as any Morita equivalent groupoid will always be a proper foliation groupoid (see
\cite[pages 18 and 21]{ademleidaruan} and \cite{crainicmoerdijk2} for more details).

Let $\Gamma$ be a finitely generated discrete group---although many
of our constructions make sense for arbitrary $\Gamma$, we are only
interested in this case. If $\phi$ and $\psi$ are homomorphisms from
$\Gamma$ to $G$, we say $\phi \sim \psi$ if they are pointwise
conjugate; i.e. if there is a $g \in G$ such that
$g\phi(\gamma)g^{-1} = \psi(\gamma)$ for each $\gamma \in \Gamma$.
We let $(\phi)$ denote the conjugacy class of $\phi$ (or sometimes
$(\phi)_\sim$ to distinguish from equivalence classes via other
relations), and let $t_{M;G}^\Gamma$ denote the set of conjugacy
classes of homomorphisms $\phi$ whose images have nonempty
fixed-point sets in $M$.  We let $M^{\langle \phi \rangle}$ denote
the fixed-point set of the image of $\phi$ in $G$ and $C_G(\phi)$ is
the centralizer of the image of $\phi$ in $G$.

\begin{definition}
\label{def-globalgammasectors}

Let $\phi : \Gamma \rightarrow G$ be a homomorphism with $M^{\langle
\phi \rangle} \neq \emptyset$. Then the {\bf $\Gamma$-sector of
$G \ltimes M$ associated to $(\phi)$} is the orbifold with presentation
\[
    (M; G)_{(\phi)} := C_G(\phi) \ltimes M^{\langle \phi \rangle}.
\]

We let $(M;G)_\Gamma$ denote the disjoint union of the
$\Gamma$-sectors,
\[
    (M;G)_\Gamma := \coprod\limits_{(\phi) \in t_{M;G}^\Gamma} (M;G)_{(\phi)} .
\]

\end{definition}

If $G$ is finite, it is obvious that each $(M; G)_{(\phi)}$ is an
orbifold groupoid (i.e. a proper \'{e}tale Lie groupoid).  We will
see in Corollary \ref{cor-sectorswelldefgoodorb} that this is
generally the case.

If $x \in M^{\langle \phi \rangle} \subseteq M$, we will sometimes
use the notation $(x, \phi)$ to distinguish between $(x, \phi) \in
M^{\langle \phi \rangle}$ and $(x, 1) \in M^{\langle 1 \rangle} =
M$.  Hence, we use $C_G(\phi)(x, \phi)$ to denote the corresponding
point in $(M;G)_{(\phi)}$

The following lemma, whose proof is standard, ensures that the definition of $(M;G)_{(\phi)}$
does not depend on the choice of the representative of the class
$(\phi)$.

\begin{lemma}
\label{lem-globgamsecswelldefined}

Let $G$ be a group acting on the smooth manifold $M$ such that $G
\ltimes M$ presents a smooth orbifold and let $\Gamma$ be a finitely
generated discrete group. If $\phi, \psi : \Gamma \rightarrow G$ are
conjugate homomorphisms with $\psi = g\phi g^{-1}$ for $g \in G$,
then the map
\[
\begin{array}{rccl}
    L_g :&  M^{\langle \phi \rangle}
        &\rTo&   M^{\langle \psi \rangle}
            \\\\
        :&  (x, \phi)  &\rMapsto&   (gx, g\phi g^{-1}) = (gx, \psi)
\end{array}
\]
is a $C_G(\phi)$-$C_G(\psi)$-equivariant diffeomorphism that induces a groupoid
isomorphism between $(M;G)_{(\phi)}$ and $(M;G)_{(\psi)}$. Moreover,
$\sigma|_{M^{\langle \phi \rangle}} = \sigma|_{M^{\langle \psi
\rangle}} \circ L_g$.

\end{lemma}

Note in particular that $G$ acts on the set $\coprod_{\phi
\in \mbox{\scriptsize HOM}(\Gamma, G)} \left(M^{\langle \phi
\rangle}, \phi\right)$ by defining $g(x, \phi) = (gx, g\phi g^{-1})$.  The
following lemma introduces a different presentation for $(M;
G)_\Gamma$.

\begin{lemma}
\label{lem-globalsectorglobalactiondef}

Suppose $G \ltimes M$ presents a quotient orbifold and let $\Gamma$ be a finitely
generated discrete group. There is a strong equivalence
\[
    G \ltimes \coprod\limits_{\psi \in \mbox{\scriptsize HOM}(\Gamma, G)}
    \left(M^{\langle \psi \rangle}, \psi\right)
    \rTo
    (M;G)_\Gamma .
\]
Hence, $G \ltimes \coprod_{\psi \in \mbox{\scriptsize
HOM}(\Gamma, G)} \left(M^{\langle \psi \rangle}, \psi\right)$ and
$(M;G)_\Gamma$ are Morita equivalent.

\end{lemma}

By strong equivalence, we mean an equivalence of groupoids such that
the map on objects is a surjective submersion; see \cite[page
20]{ademleidaruan}.  Note that neither of the groupoids in question
need be orbifold groupoids; we will see in Section \ref{sec-localvsglobaldef}
that they are both Morita equivalent to orbifold groupoids.

\begin{proof}

Pick $\phi \in \mbox{HOM}(\Gamma, G)$ with $M^{\langle \phi \rangle}
\neq \emptyset$.  Here, we denote points in $M^{\langle \phi
\rangle}$ simply as $x$ to distinguish from points in $\coprod_{\psi
\in (\phi)} \left(M^{\langle \psi \rangle}, \psi\right)$. For each
$\psi \in (\phi)$, pick a $g_\psi \in G$ such that $g_\psi \psi
g_\psi^{-1} = \phi$. We require that $g_\phi  = 1$.  Define the map
\[
\begin{array}{rccl}
    \Psi_0^\phi :& \coprod\limits_{\psi \in (\phi)} \left(M^{\langle \psi \rangle}, \psi\right)
            &\rTo&
            M^{\langle \phi \rangle}                              \\\\
        :&  (x, \psi)       &\rMapsto&       g_\psi x .
\end{array}
\]
Similarly, as $G$ acts on $\coprod_{\psi \in (\phi)} \left(M^{\langle \psi \rangle}, \psi\right)$,
define
\[
\begin{array}{rccl}
    \Psi_1^\phi :& G \times \coprod\limits_{\psi \in (\phi)} \left(M^{\langle \psi \rangle}, \psi\right)
            &\rTo&
            C_G(\phi) \times M^{\langle \phi \rangle}                              \\\\
        :&  (g, (x, \psi))       &\rMapsto&       (g_{(g\psi g^{-1})} g g_\psi^{-1}, g_\psi x).
\end{array}
\]
It is easy to check that $\Psi_0^\phi$ and $\Psi_1^\phi$ are smooth,
and that they form the maps on objects and arrows, respectively of a groupoid homomorphism
$\Psi^\phi :
G \ltimes \coprod_{\psi \in (\phi)} \left(M^{\langle \psi
\rangle}, \psi\right) \rightarrow C_G(\phi)\ltimes M^{\langle \phi
\rangle}$.

As $\Psi_0^\phi$ is a disjoint union of diffeomorphisms, $\Psi_0^\phi$ is a
surjective submersion.  It remains to show that
\[
\begin{diagram}[width=4cm]
    G \times \coprod\limits_{\psi \in (\phi)} \left(M^{\langle \psi \rangle}, \psi\right)
        &   \rTo^{\Psi_1^\phi}       &       C_G(\phi) \times M^{\langle \phi \rangle}
        \\
    \dTo<{s_1 \times s_2}            &                           &      \dTo>{s_1\times t_1}\\
    \coprod\limits_{\psi \in (\phi)} \left(M^{\langle \psi \rangle}, \psi\right)
    \times \coprod\limits_{\psi \in (\phi)} \left(M^{\langle \psi \rangle}, \psi\right)
                                    &   \rTo^{\Psi_0^\phi \times \Psi_0^\phi} &
                                    M^{\langle \phi \rangle} \times M^{\langle \phi \rangle}\\
\end{diagram}
\]
is a fibered product of manifolds.  This follows from the fact that the map
\[
\begin{array}{rccl}
    \Phi^\phi :& G \times \coprod\limits_{\psi \in (\phi)} \left(M^{\langle \psi \rangle}, \psi\right)
            &\rTo&
            \left(\coprod\limits_{\psi \in (\phi)} \left(M^{\langle \psi \rangle}, \psi\right) \right)\times
            \left(\coprod\limits_{\psi \in (\phi)} \left(M^{\langle \psi \rangle}, \psi\right) \right)\times
            \left(C_G(\phi) \times M^{\langle \phi \rangle}\right)
                                          \\\\
        :&  (h, (w, \psi))       &\rMapsto&
                            ((w, \psi), (hw, h\psi h^{-1}), (g_{(h\psi h^{-1})}hg_\psi^{-1},g_\psi w )).
\end{array}
\]
is a diffeomorphism, which is easy to verify.

With this, we need only note that $G \ltimes \coprod_{\phi
\in \mbox{\scriptsize HOM}(\Gamma, G)} \left(M^{\langle \phi
\rangle}, \phi\right)$ admits a decomposition into disjoint
groupoids
\[
    G \ltimes \coprod\limits_{\phi \in \mbox{\scriptsize HOM}(\Gamma, G)}
    \left(M^{\langle \phi \rangle}, \phi\right)
    =
    \coprod\limits_{(\phi) \in t_{M; G}^\Gamma}
    G \ltimes \coprod\limits_{\psi \in (\phi)}
    \left(M^{\langle \psi \rangle}, \psi\right),
\]
and each $\Psi^\phi$ maps one of these groupoids into $(M;
G)_{(\phi)}$.  Hence,
\[
    \Psi :=
    \coprod\limits_{(\phi) \in t_{M; G}^\Gamma} \Psi^\phi :
    G \ltimes \coprod\limits_{\phi \in \mbox{\scriptsize HOM}(\Gamma, G)}
    \left(M^{\langle \phi \rangle}, \phi\right)
    \rTo
    (M; G)_\Gamma
\]
is clearly surjective, and therefore is a strong equivalence.

\end{proof}

Note that the maps $\Psi^\phi$ depend on the choice of the $g_\psi
\in G$. However, it is easy to see that the induced map on orbit spaces does
not depend on this choice.

Fix a homomorphism $\phi : \Gamma \rightarrow G$.  Then the
injection $M^{\langle \phi \rangle} \hookrightarrow M$ induces a map
\[
\begin{array}{rccl}
    \pi_{(\phi)}    :& M^{\langle \phi \rangle}/C_G(\phi)  &\rightarrow&   M/G
    \\\\
                    :& C_G(\phi)(x,\phi)  &\mapsto&   Gx
\end{array}
\]
If $g\phi g^{-1} = \psi$ and $(x,\phi) \in M^{\langle \phi
\rangle}$, then the $G$-orbit of $x$ in $M$ coincides with that of
the corresponding point $g(x,\psi) \in M^{\langle \psi \rangle}$.
Therefore, this map does not depend on the particular choice of
representative from the conjugacy class $(\phi)$. If $(x,\phi) \in
M^{\langle \phi \rangle}$, then $\sigma(x) =
\pi_{(\phi)}\left(C_G(\phi)(x,\phi) \right)$. In particular,
$\pi_{(\phi)} \left(M^{\langle \phi \rangle} \right) = \sigma
\left(M^{\langle \phi \rangle} \right)$.

Finally, note that the map $M/C_G(\phi) \rightarrow M/G$ defined by
$C_G(\phi)(x,\phi) \mapsto Gx$ is an orbifold cover by definition
(see \cite[Definition 2.16]{ademleidaruan}). The map
$\pi_{(\phi)}$ is the restriction of this orbifold cover to
$M^{\langle \phi \rangle}/C_G(\phi)$.


\subsection{$\Gamma$-Sectors for a General Presentation}
\label{subsec-localdef}

In this subsection, we review the construction of the
$\Gamma$-sectors for a general orbifold $Q$.  We state the
construction in general for an arbitrary orbifold groupoid
$\mathcal{G}$. Throughout, we use the convention that the groupoid
$\mathcal{G}$ has space of objects $G_0$ and space of arrows $G_1$.
We also let $\sigma : G_0 \rightarrow |\mathcal{G}| = Q$ denote the
quotient map.

If $\Gamma$ and $\mathcal{G}$ are groupoids (with no additional
hypotheses), then let $\mathcal{S}_\mathcal{G}^\Gamma$ denote the
set of groupoid homomorphisms $\phi : \Gamma \rightarrow
\mathcal{G}$ such that the map on objects is constant. Then
$\mathcal{G}$ acts on $\mathcal{S}_\mathcal{G}^\Gamma$ by
conjugation; if $\phi_0 (z) = x$ for each $z \in \Gamma_0$, then for
each $g \in G_1$ with $s(g) = x$, we let $(g\cdot \phi) : \Gamma
\rightarrow \mathcal{G}$ have constant map on objects with value
$t(g)$ and map on arrows $(g\cdot \phi)_1(\gamma) =
g\phi_1(\gamma)g^{-1}$ for each $\gamma \in \Gamma_1$.

If $\Gamma$ is a group (treated as a groupoid with one unit), then
every homomorphism $\Gamma \rightarrow \mathcal{G}$ is constant on
objects and corresponds to choice of $x \in G_0$ and group
homomorphism $\phi_x :\Gamma \rightarrow G_x$ where $G_x$ denotes
the isotropy group of $x$.  Hence, we use $\phi_x$ to denote the
corresponding groupoid homomorphism.

If $\mathcal{G}$ is a topological groupoid presenting an orbispace
$X$ (see \cite{chenorbispace} or \cite{henriquesorbispace}), then each
point $x \in G_0$ is contained in an open, connected, locally connected
$U \subseteq G_0$ such that $\mathcal{G}|_U$ is isomorphic to $G_U \ltimes U$ where $G_U$ is
a topological group acting continuously on $U$.  We give
$\mathcal{S}_\mathcal{G}^\Gamma$ the weak
topology induced by the maps $\beta_\mathcal{G}: \phi_x \mapsto x \in G_0$
and for each $\gamma \in \Gamma$ the evaluation
$\epsilon_\gamma : \phi_x \mapsto \phi_x(\gamma) \in G_1$.
It is easy to check that the $\mathcal{G}$-action on $\mathcal{S}_\mathcal{G}^\Gamma$
is continuous.  With this, we make the following.

\begin{definition}
\label{def-orbispacesectors}

Let $\mathcal{G}$ be a topological groupoid representing an orbispace $X$
and let $\Gamma$ be a finitely generated discrete group.
The \emph{$\Gamma$-sector groupoid of $\mathcal{G}$}, denoted $\mathcal{G}^\Gamma$,
is the translation groupoid $\mathcal{G} \ltimes \mathcal{S}_\mathcal{G}^\Gamma$.

\end{definition}

For each $\phi_x \in \mathcal{S}_\mathcal{G}^\Gamma$ choosing
$x \in U \subseteq \mathcal{S}_\mathcal{G}^\Gamma$ as above induces an isomorphism
of topological groupoids between $C_{G_U}(\phi_x) \ltimes U^{\langle \phi_x \rangle}$
and the restriction of $\mathcal{G}^\Gamma$ to the connected component of $\beta_\Gamma^{-1}(U)$
containing $\phi_x$.
It follows that the $\Gamma$-sector groupoid represents an orbispace $|\mathcal{G}^\Gamma|$.
As a set,
\begin{equation}
\label{eq-orbispacesectors}
    |\mathcal{G}^\Gamma| = \left\{ \left(p, (\phi_x)_{\stackrel{G_x}{\sim}}
    \right) : p = \mathcal{G}x \in |\mathcal{G}|, \phi_x \in \mbox{HOM}(\Gamma, G_x) \right\}
\end{equation}
where $(\phi_x)_{\stackrel{G_x}{\sim}}$ denotes the conjugacy class of the homomorphism
$\phi_x$ in $G_x$.

Now assume that $\mathcal{G}$ is an orbifold groupoid presenting the
orbifold $Q$ and $\Gamma$ is a finitely generated discrete group.
Then if $\phi_x, \psi_y \in \mathcal{S}_\mathcal{G}^\Gamma$, a
natural transformation from $\phi_x$ to $\psi_y$ is simply a choice
of an arrow $g \in G_1$ such that $s(g) = x$, $t(g) = y$, and
$\psi_x(\gamma)g = g\phi_y(\gamma)$ for each $\gamma \in \Gamma$.
Moreover, if $\epsilon : \mathcal{K} \rightarrow \Gamma$ is an
equivalence, then $\epsilon$ is locally invertible, and $\phi_x
\circ \epsilon^{-1}$ is equivalent to $\phi_x$ (see \cite[Example
2.42]{ademleidaruan}). It follows that the orbits of points in
$\mathcal{S}_\mathcal{G}^\Gamma$ via the $\mathcal{G}$-action
correspond exactly to groupoid morphisms from $\Gamma$ to
$\mathcal{G}$.

For each point $p \in Q$ corresponding to the orbit of $x \in G_0$,
there is a linear orbifold chart $\{ V_x, G_x, \pi_x \}$ for $Q$ at
$x$.  By this, we mean that $V_x \subseteq G_0$ is diffeomorphic to
$\R^n$ with $x$ corresponding to the origin, $G_x$ acts linearly on
$V_x$, and there is a groupoid isomorphism between
$\mathcal{G}|_{V_x}$ and $G_x \ltimes V_x$.  We let $\xi_x : (s,
t)^{-1}(V_x \times V_x) \rightarrow G_x$ denote the identification
given by this isomorphism and $\xi_x^y = (\xi_x)|_{G_y} : G_y
\rightarrow G_x$ the injective homomorphism given by restriction to
$G_y$ for each $y \in V_x$.

In this case, ${\mathcal S}_{\mathcal G}^\Gamma$ is a smooth
manifold (with connected components of different dimensions) and
that the $\mathcal{G}$-action is smooth.  Hence the translation
groupoid $\mathcal{G}^\Gamma = \mathcal{G} \ltimes
\mathcal{S}_\mathcal{G}^\Gamma$ is an orbifold groupoid, defining an
orbifold structure for the {\bf $\Gamma$-sectors of $Q$}, denoted
$\tilde{Q}_\Gamma$. For each $\phi_x \in
\mathcal{S}_\mathcal{G}^\Gamma$, there is a a diffeomorphism
$\kappa_{\phi_x}$ of $V_x^{\langle \phi_x \rangle}$ onto a
neighborhood of $\phi_x$ in ${\mathcal S}_{\mathcal G}^\Gamma$
forming a manifold chart.  Identifying $V_x^{\langle \phi_x
\rangle}$ with its image via $\kappa_{\phi_x}$, $\left\{
V_x^{\langle \phi_x \rangle}, C_{G_x}(\phi_x), \pi_x^{\phi_x}
\right\}$ forms a linear orbifold chart for $\tilde{Q}_\Gamma$ at
$\phi_x$.

Within a linear chart $\{ V_x, G_x, \pi_x \}$ at $x$ with $y \in
V_x$, we say that $\phi_x$ locally covers $\psi_y$ (written $\phi_x
\stackrel{loc}{\curvearrowright} \psi_y$) if there is a $g \in G_x$
such that $g [(\xi_x^y \circ \psi_y)(\gamma)] g^{-1} =
\phi_x(\gamma)$.  Then by \cite[Lemma 2.7]{farsiseaton1}, there is a
$\psi_{y^\prime} \in \mathcal{G}\psi_y$ such that $\xi_x^y \circ
\psi_y = \phi_x$.  Extending this to an equivalence relation on
$\mathcal{S}_\mathcal{G}^\Gamma$, we say that $\phi_x \approx
\psi_y$ if there is a finite sequence of local coverings (in either
direction) connecting an element of $\mathcal{G}\phi_x$ to
$\mathcal{G}\psi_y$.  We let $(\phi)_\approx$ denote the
$\approx$-class of $\phi$ and $T_Q^\Gamma$ denote the set of
$\approx$-classes in $\mathcal{S}_\mathcal{G}^\Gamma$; when there is
no risk of confusion, we simply denote the $\approx$-class of $\phi$
by $(\phi)$. The $\approx$-classes in
$\mathcal{S}_\mathcal{G}^\Gamma$ correspond exactly to the connected
components of $\tilde{Q}_\Gamma$, so for each $(\phi) \in
T_Q^\Gamma$, we let $\tilde{Q}_{(\phi)}$ denote the connected
component consisting of $\mathcal{G}$-orbits of elements of $(\phi)$
and refer to $\tilde{Q}_{(\phi)}$ as the {\bf $\Gamma$-sector
associated to $(\phi)$}.

Note that in \cite[Lemma 2.5]{farsiseaton1}, it was shown that a
strong equivalence between orbifold groupoids induces a strong
equivalence between their associated groupoids of $\Gamma$-sectors.
Here, we will be interested in foliation groupoids that are not
necessarily \'{e}tale.  Hence, we note the following.

\begin{lemma}
\label{lem-moritaviaetale}

Suppose $\mathcal{G}$ and $\mathcal{G}^\prime$ are Morita equivalent orbifold groupoids.  Then they are
Morita equivalent via orbifold groupoids; i.e. there is an orbifold groupoid $\mathcal{H}$ and
strong equivalences
\[
    \mathcal{G}
    \lTo^e
    \mathcal{H}
    \rTo^e^\prime
    \mathcal{G}^\prime
\]

\end{lemma}

Of course, such an $\mathcal{H}$ always exists, and it is always a proper foliation groupoid.  The
point of this lemma is that $\mathcal{H}$ can taken to be \'{e}tale.

\begin{proof}

Choosing open covers of the spaces of objects consisting of linear
orbifold charts, the groupoids $\mathcal{G}$ and
$\mathcal{G}^\prime$ each give an orbifold atlas for the orbifold
$Q$ presented by $\mathcal{G}$ and $\mathcal{G}^\prime$. These
atlases need not be effective, but as they arose from a orbifold
groupoids, the kernels of the actions are appropriately restricted.
Let $\mathcal{H}$ be the groupoid of the maximal atlas containing
these two atlases, and then there are clearly equivalences as
required. Moreover, these equivalences are strong, as the domains of
charts from $\mathcal{G}$ and $\mathcal{G}^\prime$ are subsets of
the space of objects of $\mathcal{H}$ so that the embeddings of
these charts into the objects of $\mathcal{G}$ and
$\mathcal{G}^\prime$, respectively, are surjective.

\end{proof}

We solidify some notation to distinguish between the structure maps
and arrows of the groupoids under consideration.  We use $s$, $t$,
$i$, $u$, and $m$ to denote the source, target, inverse, unit, and
composition maps of a groupoid.  Often times, we will suppress $m$
and simply express products multiplicatively by concatenation; i.e.
$m(a, b) = ab$.  When it is helpful to distinguish between structure
maps of groupoids under consideration, we will give them subscripts
of the corresponding groupoid unless otherwise indicated. For a
translation groupoid $\mathcal{G}\ltimes M$, we will use the
notation throughout that $s_{\mathcal{G}\ltimes M}$ and
$t_{\mathcal{G} \ltimes M}$ are the source and target maps,
respectively, and $(\mathcal{G}\ltimes M)_1$ is the space of arrows;
note that $M$ is the space of objects. An arrow in
$(\mathcal{G}\ltimes M)_1$ is given by a $g \in G_1$ and a $z \in M$
such that the anchor map sends $s(g)$ to $z$.  We will use $(g, z)$
to denote this arrow so that $s_{\mathcal{G}\ltimes M}(g, z) = z$
and $t_{\mathcal{G}\ltimes M}(g, z) = gz$. In particular, for the
groupoid $\mathcal{G}^\Gamma =
\mathcal{G}\ltimes\mathcal{S}_\mathcal{G}^\Gamma$, an arrow is of
the form $(g, \phi_x)$ with $s_{\mathcal{G}^\Gamma}(g, \phi_x) =
\phi_x$ and $t_{\mathcal{G}^\Gamma}(g, \phi_x) = g\phi_x g^{-1}$ so
that $s(g) = x$ and $t(g) = gx$.

The following lemma will simplify many of our arguments; for the
definitions, see \cite[Definition 2.14 and 2.15]{ademleidaruan}. The
proof is direct and left to the reader.

\begin{lemma}
\label{lem-homonGspaces}

Let $\mathcal{G}$ be a groupoid, and let $M_1$ and $M_2$ be
$\mathcal{G}$-spaces with anchor maps $\alpha_i : M_i \rightarrow
G_0$. Let $e_0 : M_1 \rightarrow M_2$ be a map that is
$\mathcal{G}$-equivariant; i.e. $\alpha_2 \circ e_0 = \alpha_1$ and
$e_0(h z) = h e_0(z)$ for each $z \in M_1$ and $h \in G_1$ with
$s(h) = \alpha_1(z)$.  Define
\[
\begin{array}{rccl}
    e_1     :&      (\mathcal{G} \ltimes M_1)_1
                        &\rTo& (\mathcal{G} \ltimes M_2)_1
                                    \\\\
                   :&       (g, z)
                            &\rMapsto&   (g, e_0(z)) ,
\end{array}
\]
and then $e_0$ is the map on objects and $e_1$ the map on arrows of
a homomorphism of groupoids $e : \mathcal{G}\ltimes M_1 \rightarrow
\mathcal{G}\ltimes M_2$.  If $e_0$ is a bijection, then $e$ is an
isomorphism.

If $\mathcal{G}$ is an orbifold groupoid, the $M_i$ are smooth
$\mathcal{G}$-spaces, and $e_0$ is smooth, then $e$ is a
homomorphism of orbifold groupoids.  If $e_0$ is a diffeomorphism,
then $e$ is an isomorphism of orbifold groupoids.

\end{lemma}


\section{Connections between Definitions of Sectors}
\label{sec-localvsglobaldef}

In this section, we compare the constructions of the
$\Gamma$-sectors in Section \ref{sec-definitions} with one another,
as well as with other constructions of sectors in the literature.


\subsection{Good Orbifold}
\label{subsec-comparegood}

Let $Q$ be a good orbifold given by the quotient of a smooth manifold $M$ by a discrete group $G$
acting properly discontinuously.
Then the translation groupoid $\mathcal{G} := G \ltimes M$ is an orbifold
groupoid presenting $Q$, and $Q$ admits two
decompositions into $\Gamma$-sectors.

As in Subsection \ref{subsec-globaldef}, we let $(M;G)_\Gamma$
denote the space of $\Gamma$-sectors of $Q$ defined using the
global $G$-action on $M$; i.e. $(M;G)_\Gamma$ is given by
$\coprod_{(\phi) \in t_{M; G}^\Gamma} C_G(\phi) \ltimes
M^{\langle \phi \rangle}$.  As in Subsection \ref{subsec-localdef},
we let $\tilde{Q}_\Gamma$ denote the space of $\Gamma$-sectors of
$Q$ presented by $\mathcal{G}^\Gamma = \mathcal{G} \ltimes
\mathcal{S}_\mathcal{G}^\Gamma$.  We claim the following.

\begin{theorem}
\label{thrm-localglobaldiffeo}

Let $Q$ be a good orbifold so that $\mathcal{G} = G \ltimes M$ is an orbifold
groupoid presenting $Q$, and let $\Gamma$ be a finitely generated
discrete group.  Then $\mathcal{G}^\Gamma$ is isomorphic as an orbifold groupoid to $G \ltimes
\coprod_{\phi \in \mbox{\scriptsize HOM}(\Gamma,
G)}\left( M^{\langle \phi \rangle}, \phi \right)$.

\end{theorem}

It follows that the spaces $(M;G)_\Gamma$ and $\tilde{Q}_\Gamma$ are
diffeomorphic as orbifolds.  Before proceeding with the proof of
this proposition, we note that these spaces are not indexed in the
same way; the set $t_{M;G}^\Gamma$ is smaller than $T_Q^\Gamma$
whenever there is a homomorphism $\phi : \Gamma \rightarrow G$ such
that $\sigma(M^{\langle \phi \rangle})$ is not connected.

\begin{example}
\label{ex-localglobaldiffer}

Let $\Z/3\Z = \langle \alpha \rangle$ act on $S^2$ by rotations; the
quotient orbifold $Q$ presented by $\Z/3\Z \ltimes S^2$ is a football with two singular points,
$p_s$ and $p_n$, both of which with $\Z/3\Z$ isotropy.  Let $\Gamma =
\Z = \langle \gamma \rangle$, and define
\[
\begin{array}{rccl}
    \phi_0, \phi_1, \phi_2
                    :& \Z     &  \rTo &    \Z/3\Z       \\
    \phi_0          :& \gamma       &   \rMapsto     &   1          \\
    \phi_1          :& \gamma       &   \rMapsto     &   \alpha     \\
    \phi_2          :& \gamma       &   \rMapsto     &   \alpha^2 .
\end{array}
\]
Then the $\sim$-classes of the $\phi_i$ are the only elements of $t_{M;G}^\Gamma$.
Clearly, $(M;G)_{(\phi_0)}$ is diffeomorphic to $Q$, and
$(M;G)_{(\phi_1)}$ and $(M;G)_{(\phi_2)}$ are each diffeomorphic to
$\{ p_s, p_n \}$ with trivial $\Z/3\Z$-action.

Now, consider $\mathcal{G} \ltimes \mathcal{S}_\mathcal{G}^\Gamma$.
Let $\alpha_s$ generate $G_{p_s}$ and $\alpha_n$ generate $G_{p_n}$
for a choice of representatives of these isotropy groups.  There are
five $\approx$-classes of homomorphisms from $\Gamma$ into the local
groups of $Q$ with the following representatives:
\[
\begin{array}{rccl}
    \psi_0              :& \Z     &  \rTo &    G_p \;\;\forall p \in Q      \\
    \psi_0              :& \gamma       &   \rMapsto     &   1                  \\\\
    \psi_{1,s}         :& \Z     &  \rTo &    G_{p_s}                  \\
    \psi_{1,s}          :& \gamma       &   \rMapsto     &   \alpha_s           \\\\
    \psi_{2,s}          :& \Z     &  \rTo &    G_{p_s}                  \\
    \psi_{2,s}          :& \gamma       &   \rMapsto     &   \alpha_s^2           \\\\
    \psi_{1,n}          :& \Z     &  \rTo &    G_{p_n}                  \\
    \psi_{1,n}          :& \gamma       &   \rMapsto     &   \alpha_n           \\\\
    \psi_{2,n}          :& \Z     &  \rTo &    G_{p_n}                  \\
    \psi_{2,n}          :& \gamma       &   \rMapsto     &   \alpha_n^2 .
\end{array}
\]
Then $\tilde{Q}_{(\psi_0)}$ is diffeomorphic to $Q$, while the
sectors associated to each of the other four classes are given by a
point with trivial $\Z/3\Z$-action.

Clearly, these two decompositions result in diffeomorphic orbifolds,
although the individual sectors are indexed differently.

\end{example}

\begin{proof}[Proof of Theorem \ref{thrm-localglobaldiffeo}]

Let ${\mathcal G}$ denote the translation groupoid $G \ltimes M$ so
that $G_0 = M$ and $G_1 = G \times M$.  Then ${\mathcal G}$ is
an orbifold groupoid in the Morita equivalence class of
orbifold structures for $Q$.  We let $\zeta : G_1 = G \times M \rightarrow G$
denote the projection onto the first factor, and then
for each $\phi_x \in {\mathcal S}_{\mathcal G}^\Gamma$, we have
$\zeta \circ \phi_x \in \mbox{HOM}(\Gamma, G)$; i.e.
\[
    \zeta \circ \phi_x : \Gamma \rTo^{\phi_x} (G_1)_x
    \rTo^\zeta G.
\]
We define the map
\[
\begin{array}{rccl}
    \mathcal{Z}  :& {\mathcal S}_{\mathcal G}^\Gamma    &\rTo    &
    \coprod\limits_{\psi \in \mbox{\scriptsize HOM}(\Gamma, G)} \left( M^{\langle \psi
    \rangle}, \psi \right)
                    \\\\
            :& \phi_x               &\rMapsto
                                    &   (x,\zeta\circ\phi_x) \in \left( M^{\langle \zeta\circ\phi_x
                                    \rangle},
                                    \zeta\circ\phi_x \right)
\end{array}
\]
Then $\mathcal{Z}$ is clearly injective; if $\mathcal{Z}(\phi_x)
= \mathcal{Z}(\psi_y)$, then $(x, \zeta \circ \phi_x) = (y, \zeta
\circ \psi_y)$ so that $x = y$ and $\phi_x = \psi_y$.  To
show that $\mathcal{Z}$ is surjective, let $(x, \psi) \in
\left(M^{\langle \psi \rangle}, \psi\right)$ for some $\psi \in
\mbox{HOM}(\Gamma, G)$ and define $\psi_x : \Gamma \rightarrow
\mathcal{G}$ by $\psi_x(\gamma) = (\psi(\gamma), x)$.  Then clearly
$\mathcal{Z}(\psi_x) = (x, \psi)$, and $\mathcal{Z}$ is a bijection.
Moreover, given a chart $\kappa_{\phi_x} : V_x^{\langle \phi_x
\rangle} \rightarrow \mathcal{S}_\mathcal{G}^\Gamma$ for
$\mathcal{S}_\mathcal{G}^\Gamma$ near $\phi_x$, we have that
\[
    M^{\langle \phi_x \rangle} \supseteq
    V_x^{\langle \phi_x \rangle}
    \rTo^{\kappa_{\phi_x}}
    \mathcal{S}_\mathcal{G}^\Gamma
    \rTo^{\mathcal{Z}}
    \left( M^{\langle \zeta\circ\phi_x \rangle}, \zeta\circ\phi_x \right)
\]
is simply the identity on $V_x^{\langle \phi_x \rangle}$. It follows that
$\mathcal{Z}$ is smooth with smooth inverse, hence a diffeomorphism.

The anchor map of the $\mathcal{G}$-action on $\mathcal{S}_\mathcal{G}^\Gamma$ is
$\beta_\Gamma : \mathcal{S}_\mathcal{G}^\Gamma \rightarrow M$, with
$\beta_\Gamma : \phi_x \mapsto x$.  Let $\alpha
: \coprod_{\psi \in \mbox{\scriptsize HOM}(\Gamma, G)}
\left( M^{\langle \psi \rangle}, \psi \right) \rightarrow M$ be
defined by $\alpha: (x, \phi)\mapsto x$, and then $\alpha$ is the
anchor map of a $\mathcal{G}$-action on $\coprod_{\psi \in
\mbox{\scriptsize HOM}(\Gamma, G)} \left( M^{\langle \psi \rangle},
\psi \right) \rightarrow M$ defined by
\[
    (g, x)(x, \phi) =   (gx, g \phi g^{-1})
\]
that clearly coincides with the $G$-action.  Hence, we need only
note that for each $(g, (x, \phi)) \in G\times
\coprod_{\psi \in \mbox{\scriptsize HOM}(\Gamma, G)} \left(
M^{\langle \psi \rangle}, \psi \right)$ and $\phi_x \in
\mathcal{S}_\mathcal{G}^\Gamma$ given by $\gamma \mapsto
(\phi(\gamma),x)$ (so that $\mathcal{Z}(\phi_x) = (x, \phi)$),
\[
\begin{array}{rcl}
    (g, (x, \phi)) \mathcal{Z}(\phi_x)
        &=&     (g, (x, \phi))(x, \zeta\circ\phi_x)             \\\\
        &=&     (g, (x, \phi))(x, \phi)                         \\\\
        &=&     (gx, g\phi g^{-1})                              \\\\
        &=&     \mathcal{Z}[(g, \phi_x) \phi_x],
\end{array}
\]
and then $\mathcal{Z}$ is $\mathcal{G}$-equivariant.  It follows by
Lemma \ref{lem-homonGspaces} that $\mathcal{Z}$ is the map on
objects of an isomorphism of Lie groupoids.

\end{proof}

By Lemma \ref{lem-globalsectorglobalactiondef}, we have that
$G \ltimes \coprod_{\phi \in \mbox{\scriptsize HOM}(\Gamma,
G)} \left(M^{\langle \phi \rangle}, \phi \right)$ and $(M; G)_\Gamma$ are Morita equivalent.
Hence, by virtue of \cite[Lemma 2.5]{farsiseaton1} and Lemma \ref{lem-moritaviaetale} above,
we have the following.

\begin{corollary}
\label{cor-sectorswelldefgoodorb}

Let $Q$ be good orbifold presented by $G \ltimes M$ with $G$ discrete and let
$\Gamma$ be a finitely generated discrete group. If $\mathcal{G}$ is any orbifold groupoid presenting
$Q$, then $(M; G)_\Gamma$ and $\mathcal{G}^\Gamma$ are Morita equivalent. Hence,  the two
definitions of $\Gamma$-sectors coincide.  In particular, $(M; G)_\Gamma$ is
Morita equivalent to an orbifold groupoid.

\end{corollary}

Finally, we note that the proof of
Theorem \ref{thrm-localglobaldiffeo} generalizes readily to proper \'{e}tale orbispaces.
That is, we have the following.

\begin{theorem}
\label{thrm-localglobaldiffeoorbispace}

Let $Y$ be a $T_1$ $G$-space with $G$ discrete such that the isotropy group of each point is finite,
let $\Gamma$ be a finitely generated discrete group, and let $\mathcal{G} = G \ltimes Y$.
Then $\mathcal{G}^\Gamma$ is isomorphic as a topological groupoid to
$G \ltimes \coprod_{\phi \in \mbox{\scriptsize HOM}(\Gamma,
G)}\left( M^{\langle \phi \rangle}, \phi \right)$.

\end{theorem}

\begin{proof}

Algebraically, the proof is identical to that of Theorem \ref{thrm-localglobaldiffeo}.  Based on
the note after Definition \ref{def-orbispacesectors}, the map $\mathcal{Z}$ is clearly a
homeomorphism.  Moreover, the induced map on arrows given by Lemma \ref{lem-homonGspaces}
is clearly a homeomorphism as well.

\end{proof}


\subsection{Quotient Orbifolds}
\label{subsec-comparegeneral}

In the case that $G$ is not discrete, we have the following.

\begin{theorem}
\label{thrm-localglobaldiffeogenerallie}

Let $G$ be a Lie group that acts smoothly on the smooth manifold $M$
satisfying conditions (i), (ii), and (iii) in Subsection
\ref{subsec-globaldef} so that $G \ltimes M$ presents an orbifold $Q$. Let
$\mathcal{G}$ be an orbifold groupoid representing $Q$ so that $G
\ltimes M$ and $\mathcal{G}$ are Morita equivalent. Then
$\mathcal{G}^\Gamma$ and $(M;G)_\Gamma$ are
Morita equivalent.

\end{theorem}

\begin{proof}

First, we construct a specific orbifold groupoid that is Morita
equivalent to $G \ltimes M$.

If $G$ acts properly on $M$ with discrete isotropy groups, then $M$
is foliated by (connected components of) $G$-orbits (see \cite[page
16]{moerdijkmrcun}). Pick $x \in M$, and then there is a unique
$G_x$-space $S_x$ and a $G$-diffeomorphism of $G \times_{G_x} S_x$
onto an open subset of $M$ containing $x$. We recall the
construction of $G \times_{G_x} S_x$.  If $(u, y) \in G \times S_x$
and $k \in G_x$, then $k(u, y) = (uk^{-1}, ky)$ defines a
$G_x$-action on $G \times S_x$, and $G \times_{G_x} S_x$ is the
orbit space of this action.  Then the $G$-action on $G \times S_x$
given by $g^\prime(g, y) = (g^\prime g, y)$ induces a $G$-action on
$G \times_{G_x} S_x$ (see \cite[page 32]{tomdieck}). In particular,
the slice $S_x$ is a transversal for the foliation of $(G
\times_{G_x} S_x)$ by $G$-orbits.  We note that $S_x$ is not a
complete transversal unless $G/G_x$ is connected; in general, a
complete transversal to the foliation of $(G \times_{G_x} S_x)$ can
be formed by picking one translate $gS_x$ of the slice of $S_x$ in
each connected component of $(G \times_{G_x} S_x)$.

As $M/G$ is paracompact, an open cover of $M/G$ formed by picking a
chart of the form $G \times_{G_x} S_x$ for a choice of one point $x$
in each $G$-orbit of $M$ can be refined to a locally finite cover by
shrinking the $S_x$; hence, we can form a complete transversal $S$
to the foliation of $M$ by $G$-orbits by taking the (possibly
disconnected) union of slices $S_x$.

By \cite[Theorem 1 and Lemma 2]{crainicmoerdijk2}, $G \ltimes M$ is equivalent to
the groupoid given by the restriction $(G \ltimes M)|_S$ of $G \ltimes M$ to a complete transversal $S$
(note that the \emph{essential equivalence} of \cite{crainicmoerdijk2} corresponds to an
\emph{equivalence} in \cite[Definition 1.42]{ademleidaruan}; we use the language of the latter
for consistency).  Moreover, $(G \ltimes M)|_S$ is \'{e}tale.
Since $G\ltimes M$ is a proper and properness is preserved under equivalence,
$(G \ltimes M)|_S$ is an orbifold groupoid.

The following argument follows \cite[Theorem 5.3]{ademruan}, which treats the case of $\Gamma = \Z$.

Pick a homomorphism $\phi : \Gamma \rightarrow G$ with nonempty fixed-point set in $M$.  As $G$ acts on
$M$ with discrete isotropy, $C_G(\phi)$ clearly acts on $M^{\langle \phi \rangle}$ with discrete isotropy
and hence foliates $M^{\langle \phi \rangle}$ by (connected components of) $C_G(\phi)$-orbits.  We
construct a complete transversal to this foliation from the complete transversal $S$.

Pick a chart of the form $(G \times_{G_x} S_x)$ where the slice
$S_x$ is contained in $S$.  Then $(G \times_{G_x} S_x)^{\langle \phi
\rangle}$ is by definition the set of $G_x(u, y) \in (G \times_{G_x}
S_x)$ such that
\[
    \forall \gamma \in \Gamma \;\; \exists h \in G_x :
    (\phi(\gamma)u,y) = h(u, y)
\]
where again $h(u,y) = (uh^{-1}, hy)$.  We claim that $(G \times_{G_x} S_x)^{\langle \phi \rangle}$
is given by
\begin{equation}
\label{eq-slicefixedpoints}
    \left\{ G_x(u ,y) \in (G \times_{G_x} S_x) : u^{-1} (\Imt \phi) u \leq G_x, y
    \in S_x^{\langle u^{-1}\phi u \rangle} \right\}.
\end{equation}

Suppose $u^{-1} (\Imt\phi)u \leq G_x$ and $y \in S_x^{\langle
u^{-1}\phi u \rangle}$.  For each $\gamma \in \Gamma$,
\[
\begin{array}{rcl}
    u^{-1} \phi(\gamma)^{-1} u(u, y)
        &=&     (u(u^{-1} \phi(\gamma)^{-1}u)^{-1}, u^{-1} \phi(\gamma)^{-1}uy)      \\\\
        &=&     (uu^{-1} \phi(\gamma)u, y)        \\\\
        &=&     (\phi(\gamma)u, y).
\end{array}
\]
As $u^{-1} \phi(\gamma)^{-1}u \in G_x$, it follows that the
$G_x$-orbits $G_x(\phi(\gamma)u, y) = G_x(u, y)$.  As this is true
for each $\gamma \in \Gamma$, we have that $G_x(u, y) \in (S_x
\times_{G_x} G)^{\langle \phi \rangle}$.

Conversely, suppose the orbit $G_x(u, y)$ is fixed by $\phi(\gamma)$
for each $\gamma \in \Gamma$. Then for each $\gamma \in \Gamma$,
there is an $h \in G_x$ such that $(\phi(\gamma)u, y) = h(u, y) =
(uh^{-1}, hy)$. It follows that $\phi(\gamma)u = uh^{-1}$ so that
$u^{-1} \phi(\gamma) u = h^{-1} \in G_x$. Moreover, $y = hy$ so that
$y \in S_x^{\langle h \rangle} = S_x^{\langle u^{-1}
\phi(\gamma)^{-1} u \rangle}$. As this is true for each $\gamma \in
\Gamma$, $u^{-1} (\Imt \phi) u \leq G_x$ and $y \in S_x^{\langle
u^{-1}\phi u \rangle}$, proving the expression in (\ref{eq-slicefixedpoints}) of
$(G \times_{G_x} S_x)^{\langle \phi \rangle}$.

Now, let $\mathcal{O}_\phi$ be the collection of $\psi : \Gamma
\rightarrow G_x \leq G$ that are conjugate to $\phi$ in $G$.  Then
$G_x$ acts on $\coprod_{\psi \in \mathcal{O}_\phi}
\left(S_x^{\langle \psi \rangle}, \psi \right)$ via $h(y, \psi) =
(hy, h\psi h^{-1})$.  We let $[y, \psi]$ denote the $G_x$-orbit of
$(y, \psi)$.  Define the map
\[
\begin{array}{rccl}
    \mathcal{E}:    &   (G \times_{G_x} S_x)^{\langle \phi \rangle}
                    &\rTo&
                    \left(\coprod\limits_{\psi \in \mathcal{O}_\phi}
                        \left(S_x^{\langle \psi \rangle}, \psi \right)\right)/G_x
                                \\\\
                    :&  G_x(u, y)
                    &\rMapsto&
                        [y, u^{-1}\phi u].
\end{array}
\]
This map is well-defined, as for $h \in G_x$,
\[
\begin{array}{rcl}
    \mathcal{E}(G_x h(u, y))
        &=&     \mathcal{E}(G_x(uh^{-1}, hy))                \\\\
        &=&     [hy, hu^{-1}\phi uh^{-1}]
                                      \\\\
        &=&     h[y, u^{-1}\phi u]                                      \\\\
        &=&     [y, u^{-1}\phi u]                                      \\\\
        &=&     \mathcal{E}(u, y).
\end{array}
\]
Note that $y \in S_x^{\langle u^{-1}\phi u \rangle}$ whenever
$G_x(u, y) \in (G \times_{G_x}S_x)^{\langle \phi \rangle}$, and note
further that the map $\mathcal{E}$ is clearly smooth, both
observations by virtue of (\ref{eq-slicefixedpoints}).

The map $\mathcal{E}$ is not injective.  However, we claim that
$\mathcal{E}(G_x(u, y)) = \mathcal{E}(G_x(v, y^\prime))$
if and only if there is a $z \in C_G(\phi)$ such that $z(G_x(u, y)) = G_x(v, y^\prime)$; i.e.
$(zu, y) = h(v, y^\prime) = (vh^{-1}, hy^\prime)$ for some $h \in G_x$.  If this is
the case, then $(u, y) = (z^{-1}vh^{-1}, hy^\prime)$, so that
\[
\begin{array}{rcl}
    \mathcal{E}(G_x(u,y))
        &=&     [y, u^{-1}\phi u]
                        \\\\
        &=&     [hy^\prime, (z^{-1}vh^{-1})^{-1}\phi z^{-1}vh^{-1}]
                        \\\\
        &=&     [hy^\prime, hv^{-1}z \phi z^{-1}vh^{-1}]
                        \\\\
        &=&     h[y^\prime, v^{-1} \phi v]
                        \\\\
        &=&     [y^\prime, v^{-1} \phi v]
                        \\\\
        &=&     \mathcal{E}(G_x(v, y^\prime)).
\end{array}
\]
Conversely, if $\mathcal{E}(G_x(u,y)) = \mathcal{E}(G_x(v, y^\prime))$, then
$[y, u^{-1}\phi u] = [y^\prime, v^{-1}\phi v]$ so that there is an $h \in G_x$ such that
$(y, u^{-1}\phi u) = h(y^\prime, v^{-1}\phi v) = (hy^\prime, hv^{-1}\phi vh^{-1})$.
It follows that $y = hy^\prime$ and $u^{-1}\phi u = hv^{-1} \phi vh^{-1}$; i.e. that
$\phi = vh^{-1}u^{-1}\phi uhv^{-1}$.  Hence letting $z = vh^{-1}u^{-1}$, $z \in C_G(\phi)$,
and we have that $zu = vh^{-1}$, so that $(zu, y) = (vh^{-1}, hy^\prime)$.

To see that this map is surjective, let $[y, \psi] \in
\coprod_{\psi \in \mathcal{O}_\phi} \left(S_x^{\langle \psi
\rangle}, \psi \right)$ and then there is a $u \in G$ such that
$u\psi u^{-1} = \phi$.  Then $(u, y) \in (G\times_{G_x}S_x)^{\langle
\phi \rangle}$ and $\mathcal{E}(G_x(u, y)) = [y, \psi]$.

With this, we have that $\mathcal{E}$ induces a diffeomorphism from
$(G \times_{G_x} S_x)^{\langle \phi \rangle} /C_G(\phi)$ onto
$\left(\coprod_{\psi \in \mathcal{O}_\phi} (S_x^{\langle
\psi \rangle}, \psi)\right)/G_x$. Let
$(\psi)_{\stackrel{G_x}{\sim}}$ denote the $G_x$-conjugacy class of
$\psi$ to distinguish it from the $G$-conjugacy class.  Recall from
the proof of Lemma \ref{lem-globalsectorglobalactiondef} that the
strong equivalence
\[
    G_x \ltimes \coprod\limits_{\psi \in \mbox{\scriptsize HOM}(\Gamma, G_x)} (S_x^{\langle \psi \rangle}, \psi)
    \rTo
    \coprod\limits_{(\psi)_{\stackrel{G_x}{\sim}} \in \mathcal{O}_\phi/G_x} C_{G_x}(\psi) \ltimes
    S_x^{\langle \psi \rangle}
\]
restricts to an equivalence
\[
    G_x \ltimes \coprod\limits_{\psi_0 \in (\psi)_{\stackrel{G_x}{\sim}}} (S_x^{\langle \psi \rangle}, \psi)
    \rTo
    C_{G_x}(\psi) \ltimes S_x^{\langle \psi \rangle}
\]
for each $\stackrel{G_x}{\sim}$-class $(\psi)_{\stackrel{G_x}{\sim}}$.  Noting that $\mathcal{O}_\phi$
clearly consists of entire $\stackrel{G_x}{\sim}$-classes, we have that there is an equivalence
\[
    G_x \ltimes \coprod\limits_{\psi \in \mathcal{O}_\phi} (S_x^{\langle \psi \rangle}, \psi)
    \rTo
    \coprod\limits_{(\psi)_{\stackrel{G_x}{\sim}} \in \mathcal{O}_\phi/G_x} C_{G_x}(\psi) \ltimes
    S_x^{\langle \psi \rangle},
\]
where the $G_x$-action on $\mathcal{O}_{\phi}$ is by conjugation.  This implies that there is a diffeomorphism
\begin{equation}
\label{eq-localunionofsectorscentralizeraction}
    \left(\coprod\limits_{\psi \in \mathcal{O}_\phi} (S_x^{\langle \psi \rangle}, \psi)\right)/G_x
    \rTo
    \coprod\limits_{(\psi)_{\stackrel{G_x}{\sim}} \in \mathcal{O}_\phi/G_x}
    S_x^{\langle \psi \rangle}/C_{G_x}(\psi).
\end{equation}

Note that $(G \times_{G_x} S_x)^{\langle \phi \rangle}$ is empty
unless $\phi$ is conjugate in $G$ to a homomorphism with image in
$G_x$.  Choose one representative $\psi$ from each $G_x$-conjugacy
class $(\psi)_{\stackrel{G_x}{\sim}}$.  Recall that the map
$\mathcal{E}$ is constant on $C_G(\phi)$-orbits.   From its
definition, $\mathcal{E}$ maps the submanifold $S_x^{\langle \psi
\rangle}$ of the slice $S_x$ to the $G_x$-orbit of $\left(
S_x^{\langle \psi \rangle}, \psi\right)$. Moreover, if $\psi$ is the
chosen representative of the conjugacy class
$(\psi)_{\stackrel{G_x}{\sim}}$, the equivalence in Lemma
\ref{lem-globalsectorglobalactiondef} maps $\left( S_x^{\langle \psi
\rangle}, \psi\right)$ onto $S_x^{\langle\psi\rangle}$. It follows
from this diffeomorphism and these observations that the disjoint
union $\mathcal{S}_x =
\coprod_{(\psi)_{\stackrel{G_x}{\sim}} \in
\mathcal{O}_\phi/G_x} S_x^{\langle \psi \rangle}$ is a complete
transversal to the foliation of $(G \times_{G_x} S_x)^{\langle \phi
\rangle}$ by connected components of $C_G(\phi)$-orbits.
Forming $\mathcal{S}_x$ for each chart for $S$ as above, the
(possibly disconnected) union $\tilde{S}$ of the $\mathcal{S}_x$
forms a complete transversal to the foliation of $M^{\langle \phi
\rangle}$ by the $C_G(\phi)$-action.

As usual, let $(G\ltimes M)|_S^\Gamma$ denote the groupoid of $\Gamma$-sectors for
the orbifold groupoid $(G\ltimes M)|_S$, constructed as in Subsection \ref{subsec-localdef}.
Note that the space of objects of $(G\ltimes M)|_S$ is simply $S$ while the arrows of $(G\ltimes M)|_S$
are given by $(g, x) \in G \times S$ such that $gx \in S$.  Clearly, then, the isotropy group of
a point $x \in S$ is simply $G_x$, the isotropy group of $x$ as a point in $M$.
It follows that the space of objects of $(G\ltimes M)|_S^\Gamma$ is the set of homomorphisms
$\psi_x : \Gamma \rightarrow G_x$ for $x \in S$ with local charts given by
$V_x^{\langle \psi_x \rangle}$.  As the action of an arrow $(g, x)$ in
$(G\ltimes M)|_S$ is given by $g\phi_x g^{-1}$, yielding a homomorphism from $\Gamma$ into
$G_{gx}$, the groupoid $(G\ltimes M)|_S^\Gamma$ is isomorphic to the restriction of the
groupoid $C_G(\phi) \ltimes M^{\langle \phi \rangle}$ to the complete transversal $\tilde{S}$ given
above.  As this is true for each $(\phi) \in t_{M;G}^\Gamma$, it follows that there is an equivalence from
$(G\ltimes M)|_S^\Gamma$ to $G \ltimes \coprod_{\phi \in
\mbox{\scriptsize HOM}(\Gamma, G)}M^{\langle \phi \rangle}$.  Hence
$(G\ltimes M)|_S^\Gamma$ and $(M;G)_\Gamma$ are Morita equivalent by
Lemma \ref{lem-globalsectorglobalactiondef}.

To complete the proof, suppose that $\mathcal{G}$ is any orbifold groupoid Morita equivalent to
$G\ltimes M$.  Then $\mathcal{G}$ is Morita equivalent to $(G\ltimes M)|_S$ via \'{e}tale groupoids
by Lemma \ref{lem-moritaviaetale}, implying by \cite[Lemma 2.5]{farsiseaton1} that
the $\Gamma$-sectors of the two groupoids are Morita equivalent.

\end{proof}

\begin{corollary}
\label{cor-liesectorsorbifoldsgoodorb}

Let $G$ be a Lie group that acts smoothly on the smooth manifold $M$
satisfying conditions (i), (ii), and (iii) in Subsection
\ref{subsec-globaldef} so that $G \ltimes M$ presents an orbifold $Q$.
Let $\Gamma$ be a finitely generated discrete group.
Then $(M;G)_\Gamma$ is Morita equivalent to an orbifold groupoid and hence
presents an orbifold.

\end{corollary}

While Example
\ref{ex-localglobaldiffer} illustrates that the correspondence
\[
    T_Q^\Gamma \ni (\phi_x)_\approx \rMapsto (\zeta \circ \phi_x)_\sim \in
    t_{M;G}^\Gamma
\]
is not injective, it is clearly surjective.  It is an obvious
consequence of Theorems \ref{thrm-localglobaldiffeo} and
\ref{thrm-localglobaldiffeogenerallie} and the fact that
$\approx$-classes are precisely connected components of
$\tilde{Q}_\Gamma = |\mathcal{G}^\Gamma|$ that each $\approx$-class
corresponds to a connected component of a $\sim$-class of
$(M;G)_\Gamma$.

With this, we note that the equivalence $\approx$ defined on objects
of $\mathcal{G}^\Gamma$ in Subsection \ref{subsec-localdef} can be
expressed naturally on either model of $(M;G)_\Gamma$.  Using the
groupoid $(M; G)_\Gamma$ defined in Definition
\ref{def-globalgammasectors}, given $(x, \phi), (y, \phi) \in
M^{\langle \phi \rangle}$, we say that $(x, \phi) \approx (y, \psi)$
if the orbits $C_G(\phi)x$ and $C_G(\phi)y$ are on the same
connected component of $M^{\langle \phi \rangle}/C_G(\phi)$.
Similarly, using the Morita equivalent groupoid representing
$(M;G)_\Gamma$ given by Lemma \ref{lem-globalsectorglobalactiondef},
we say that $(x, \phi) \approx (y, \psi)$ for two points $(x, \phi),
(y, \psi) \in \coprod_{\psi \in \mbox{\scriptsize
HOM}(\Gamma, G)} \left( M^{\langle \psi\rangle}, \psi\right)$
whenever there is a $g \in G$ such that $g\phi g^{-1} = \psi$ and
such that the orbits $G(x, \phi)$ and $G(y, \psi) = G(gy, \phi)$ are
on the same connected component of $\left(\coprod_{\psi \in
\mbox{\scriptsize HOM}(\Gamma, G)} \left( M^{\langle \psi\rangle},
\psi\right)\right)/G$.  Clearly, the three definitions of $\approx$
coincide in the sense that they define the same equivalence classes
on the quotient space, and the $\approx$-classes correspond exactly
to connected components.  We let $(x, \phi)_\approx$ denote the
$\approx$-class of the point $(x, \phi)$ in either case and $T_{M;
G}^\Gamma$ the set of $\approx$-classes.  Then $T_{M;G}^\Gamma$ and
$T_Q^\Gamma$ obviously coincide.

In the same way, the definitions in \cite[Section 3]{farsiseaton1}
can be reformulated from the perspective of a presentation as a
global quotient.  Let $(x, \phi)_\approx, (y, \psi)_\approx \in
T_{M;G}^\Gamma$ and let $(M^{\langle \phi \rangle}/C_G(\phi))_x$ and
$(M^{\langle\psi\rangle}/C_G(\psi))_y$ denote the connected
components of $M^{\langle \phi \rangle}/C_G(\phi)$ and
$M^{\langle\psi\rangle}/C_G(\psi)$ containing the orbits of $x$ and
$y$, respectively.  We say that $(x, \phi)_\approx \leq (y,
\psi)_\approx$ if $\pi((M^{\langle \phi\rangle}/C_G(\phi))_x) \subseteq
\pi((M^{\langle\psi\rangle}/C_G(\phi))_y)$ where $\pi : (M;G)_\Gamma \rightarrow
M/G$ denotes the map $C_G(\phi)(x, \phi) \mapsto Gx$.  Similarly,
$\Gamma$ \emph{covers the local groups of $Q$} if, for every $H \leq G$
such that $M^H \neq \emptyset$, there is a surjective homomorphism
$\phi : \Gamma \rightarrow H$.


\subsection{Connections between $\Gamma$-Sectors and Other Sectors}
\label{subsec-otherconstructions}

The definition of the $\Gamma$-sectors was motivated by that of the inertia orbifold and
the $k$-multi-sectors given in \cite[pages 52--53]{ademleidaruan} (see also \cite{chenruanorbcohom}).
Hence, the $\Gamma$-sectors generalize the definition of the multi-sectors in the following sense.

\begin{proposition}
\label{prop-multisectorsgammasectors}

Let $Q$ be an orbifold presented by the orbifold groupoid
$\mathcal{G}$ and let $\mathbb{F}_k$ denote the free group
on $l$-generators. Then the groupoids $\mathcal{G} \ltimes
\mathcal{S}_\mathcal{G}^k$ and $\mathcal{G} \ltimes
\mathcal{S}_\mathcal{G}^{\F_k}$ are isomorphic.  In particular,
$\tilde{Q}_{\mathbb{F}_k}$ is diffeomorphic to the space of
$k$-multi-sectors $\tilde{Q}_k$.

\end{proposition}

\begin{proof}

This follows almost immediately from the definition.
Let $\mathbb{F}_k$ be generated by $\gamma_1, \ldots , \gamma_k$, and
recall from \cite{ademleidaruan} that ${\mathcal S}_{\mathcal G}^k$
is defined to be the set
\[
    \{ (g_1, \ldots, g_k): g_i \in G_1, s(g_i) = t(g_j) \;\;\;\forall i, j \leq k \}.
\]
To each $(g_1, \ldots, g_k) \in {\mathcal S}_{\mathcal
G}^k$ with $s(g_i) = t(g_j) = x$, there is a unique homomorphism
$\phi_x : \mathbb{F}_k \rightarrow G_x$ such that $\phi_x(\gamma_i) =
g_i$.  It is obvious that the identification $(g_1, \ldots, g_k) \mapsto \phi_x$
is a homeomorphism ${\mathcal S}_{\mathcal G}^k \rightarrow {\mathcal
S}_{\mathcal G}^{\mathbb{F}_k}$.  With
this, we need only note that the action of ${\mathcal G}$ on
${\mathcal S}_{\mathcal G}^k$ and ${\mathcal S}_{\mathcal
G}^{\mathbb{F}_k}$ are defined identically, and hence the result
follows by an application of Lemma \ref{lem-homonGspaces}.

\end{proof}

\begin{corollary}
\label{cor-inertiaZsectors}

Let $\mathcal{G}$ be an orbifold groupoid. Then $\mathcal{G}^\Z$ is isomorphic as a
groupoid to the inertia groupoid $\wedge \mathcal{G}$.  In
particular, the space of $\Z$-sectors $\tilde{Q}_\Z$ is
diffeomorphic to the inertia orbifold $\tilde{Q}$.

\end{corollary}

In \cite{leida}, Leida defines the \emph{fixed-point sectors} of an
orbifold groupoid $\mathcal{G}$.  Recall that Leida defines
$\tilde{\mathcal{S}}(\mathcal{G}) = \{ (x, H)| x \in G_0, H \leq G_x
\}$, and $\tilde{\mathcal{G}} = \mathcal{G} \ltimes
\tilde{\mathcal{S}}(\mathcal{G})$.  Similarly, for each subgroup $H$
of $G_1$, $\tilde{\mathcal{S}}^H(\mathcal{G})$ is the subset $\{ (x,
K)| K \cong H \}$ given by points $(x, K)$ where $K$ is isomorphic
to $H$. Define the map
\[
\begin{array}{rccl}
    \varrho         :&    \mathcal{S}_\mathcal{G}^\Gamma  &\rTo&
                            \tilde{\mathcal{S}}(\mathcal{G})
                                \\\\
                    :&    \phi_x    &\rMapsto&       (x, \Imt
                    \phi_x).
\end{array}
\]
For each point $(x, \Imt \phi_x) = \varrho(\phi_x)$ in the image of
$\varrho$, there is a neighborhood $V_x$ of $x$ in $G_0$ such that
the restriction $\mathcal{G}|_{V_x}$ is isomorphic to $G_x \ltimes
V_x$. This corresponds to a neighborhood of $(x, \Imt\phi_x)$ in
$\tilde{\mathcal{S}}^{\mbox{\scriptsize Im\,} \phi_x}(\mathcal{G})$
diffeomorphic to $V_x^{\langle \phi_x \rangle}$ such that the
restriction of $\tilde{\mathcal{G}}^{\mbox{\scriptsize Im\,}
\phi_x}$ is isomorphic to $N_{G_x}(\Imt\phi_x) \ltimes V_x^{\langle
\phi_x \rangle}$ (see \cite[Section 2.2]{leida};
$N_{G_x}(\Imt\phi_x)$ denotes the normalizer of $\Imt\phi_x$ in
$G_x$). Similarly, there is a neighborhood of $\phi_x$ in
$\mathcal{S}_\mathcal{G}^\Gamma$ such that the restriction of
$\mathcal{G}^\Gamma$ is isomorphic to $C_{G_x}(\phi_x) \ltimes
V_x^{\langle \phi_x \rangle}$.  When restricted to these
neighborhoods, the map $\varrho$ is simply the embedding of
$C_{G_x}(\phi_x) \ltimes V_x^{\langle\phi_x\rangle}$ into
$N_{G_x}(\phi_x) \ltimes V_x^{\langle\phi_x\rangle}$.  If $\psi_x$
is another point with $\varrho(\psi_x) = (x, \Imt \phi_x)$ then
$\Imt \psi_x = \Imt \phi_x$ so that $C_{G_x}(\psi_x) \ltimes
V_x^{\langle\psi_x\rangle} = C_{G_x}(\phi_x) \ltimes
V_x^{\langle\phi_x\rangle}$.

If a point $(x, \Imt \phi_x)$ is in the image of $\varrho$, then
every point in $\tilde{\mathcal{G}}^{\mbox{\scriptsize Im\,}
\phi_x}$ is in the image of $\varrho$.  To see this, note that if
$(y, H) \in \tilde{\mathcal{S}}^{\mbox{\scriptsize Im\,}
\phi_x}(\mathcal{G})$, then as $H$ is isomorphic to $\Imt \phi_x$,
there is a homomorphism $\psi_y : \Gamma \rightarrow G_y$ with image
$H$.  It follows that $\varrho(\psi_y) = (y, H)$.  Note that it need
not be the case that $\psi_y \approx \phi_x$.  However, using the
techniques of \cite[Lemma 3.2]{farsiseaton1}, it is easy to see that
the images of $\Gamma$-sectors $(\phi)$ via $\varrho$ are entire
connected components of $\tilde{\mathcal{G}}$. In particular, if
$(x, \Imt \phi_x)$ and $(y, H)$ are in the same connected component
of a fixed-point sector $\tilde{\mathcal{G}}^{\mbox{\scriptsize
Im\,} \phi_x} = \tilde{\mathcal{G}}^H$, then they are connected by a
path in $\tilde{\mathcal{G}}^H$ and hence a finite number of charts of the form
$N_{G_{x_i}}(H_i) \ltimes V_{x_i}^{H_i}$ with each $H_i$ isomorphic.  Arrows $g \in G_1$
connecting orbits of point in $V_{x_i}^{H_i}$ to those in
$V_{x_{i+1}}^{H_{i+1}}$ act on homomorphisms with images in $H_i$
resulting in homomorphisms with images in $H_{i+1}$.  Hence, a
sequence of homomorphisms $\phi_{x_i}$ can be defined in each chart,
showing that there is a $\phi_y \approx \phi_x$ with
$\varrho(\phi_y) = (y, H)$.

If $\Gamma$ covers the local groups of $Q$, i.e. for each subgroup $H \leq G_x$
of an isotropy group of $Q$, there is a surjective homomorphism $\Gamma\rightarrow H$,
then it is clear that each $(x,
H)$ is the image via $\varrho$ of a $\phi_x$ with $\Imt\phi_x = H$.
Hence, $\varrho$ is surjective, and we have the following.

\begin{proposition}
\label{prop-leidasectors}

Let $\Gamma$ be a finitely generated discrete group and $\mathcal{G}$ an
orbifold groupoid presenting the orbifold $Q$.  Each $\Gamma$-sector
of $Q$ is an orbifold cover of a connected component of each
$\tilde{\mathcal{G}}^{\mbox{\scriptsize Im\,} \phi_x}(\mathcal{G})$.
If $\phi_x$ is chosen to have minimal isotropy in $(\phi)$, then the
$\Gamma$-sector $(\phi)$ is a
$[N_{G_x}(\phi_x):C_{G_x}(\phi_x)]$-cover of the corresponding
fixed-point sector.  If $\Gamma$ covers the local groups of $Q$,
then each connected component of each $\tilde{\mathcal{G}}^H$ is
orbifold covered by a $\Gamma$-sector.

\end{proposition}

As the homotopy groups of an orbifold groupoid $\mathcal{G}$ are Morita invariant,
it follows that the homotopy groups of the $\Gamma$-sectors are Morita invariants
for each finitely generated discrete group $\Gamma$.  See \cite{leida}, \cite{chenorbispace}, and
\cite{henriquesorbispace} for more on homotopy theory of homotopy groups of orbifolds and orbispaces.


\section{A Model of the $\Gamma$-Sectors Using Generalized Loop Spaces}
\label{sec-multiloop}

In \cite{luperciouribeloop}, it is shown that the inertia orbifold
of an orbifold $Q$ occurs in the context of the loop space of $Q$.
It appears as the subset of constant
loops or equivalently the set of loops fixed by the natural $S^1$-action on
the loop space.  In this section,
we show how this construction can be generalized to demonstrate that the
$\Gamma$-sectors of an orbifold arise in the same way when
considering maps from a closed manifold $M_\Gamma$ with fundamental
group $\Gamma$.  See also \cite{tamanoi2} for similar results for global quotients from
a different perspective.

Many of the results in this section can be proven by direct generalizations of
arguments in \cite{luperciouribeloop} once the appropriate definitions are given.
Hence, we will be thorough with the details of the definitions and refer the
reader to \cite{luperciouribeloop}, noting any nontrivial changes.
Throughout this section, we let $Q$ be an arbitrary smooth orbifold
represented by the orbifold groupoid $\mathcal{G}$.


\subsection{The $M_\Gamma$-Multiloop Space of an Orbifold}
\label{subsec-intromultiloop}

In this subsection, we develop a groupoid structure for a manifold
$M_\Gamma$ with fundamental group $\Gamma$.  This construction generalizes
that of \cite[Sections 3.1--3.2]{luperciouribeloop} for the case of $\Gamma = \Z$ and $M_\Gamma = S^1$.

Let $\Gamma$ be a finitely generated discrete group, $M_\Gamma$ a smooth manifold
with fundamental group $\Gamma$, and $M$ the universal cover of
$M_\Gamma$ so that $M/\Gamma = M_\Gamma$.  We let $\pi_\Gamma  : {M}
\to M_\Gamma $ denote the covering projection.
Fix a metric on $M_\Gamma$ and consider a cover ${\mathcal U} = \{
U_n \}_{n\in \mathbb N} $ of $M_\Gamma$ that is
$\frac{1}{n}$-admissible; i.e. each $U_n$ is evenly covered and has
diameter $\leq \frac{1}{n}$.  Note that if $M_\Gamma$ is compact, we
can assume that ${\mathcal U}$ is finite. Let ${\mathcal W}$ be the
cover of $M$ formed by the connected components of the sets
$\pi_\Gamma^{-1}(U_i)$ for each $U_i \in {\mathcal U}$. In other
words, for each $n \in \mathbb{N}$, choose one connected component
$W_n^1$ of $\pi_\Gamma^{-1}(U_n)$ and let $W_n^\gamma = \gamma
W_n^1$.  Then ${\mathcal W} = \{ W_n^\gamma \}_{n \in \mathbb{N}, \gamma
\in \Gamma}$.  Set $W_n = \pi^{-1}_\Gamma (U_n ) =
\coprod_{\gamma\in \Gamma} W_n^\gamma$, and define the
groupoid $M^{\mathcal W} $ to be the groupoid associated to the
covering ${\mathcal W}$ of $M$.  That is, the set of the units
$(M^{\mathcal W} )_0$ of $M^{\mathcal W} $ is given by
\[
    (M^{\mathcal W} )_0   = \coprod\limits_{n \in {\mathbb N}, \gamma\in \Gamma}
    W_n^\gamma,
\]
and the set of arrows $(M^{\mathcal W} )_1$ is
\[
    (M^{\mathcal W} )_1  = \coprod\limits_{n, m \in {\mathbb N}; \gamma, \delta
    \in \Gamma }  W_n^\gamma \cap  W_m^\delta.
\]
We let $(x, W_n^\gamma)$ denote the object associated to $x \in W_n^\gamma \subseteq M$
to distinguish it from $(x, W_m^\delta)$ in the case that $x \in W_n^\gamma \cap W_m^\delta$.
When the specific translate of $W_n^1$ does not concern us, we simply use $(x, W_n)$.
Note that this introduces no ambiguity; $\Gamma$ is the group of deck translations of the
manifold cover $M \rightarrow M_\Gamma$ so that $x$ can be contained in only one translate $W_n^\gamma$
of $W_n^1$. Similarly, we use $W_{n, m}^{\gamma, \delta}$ to denote the connected component
$W_n^\gamma \cap  W_m^\delta$ of $(M^\mathcal{W})_1$ and let
$W_{n,m}=\coprod_{\gamma,\delta\in \Gamma} W_{n,m}^{\gamma,\delta}$.
Then $(x, W_{n, m}^{\gamma, \delta})$ or simply $(x, W_{n,m})$ (again, with no ambiguity)
denotes the arrow corresponding to the point $x \in W_{n, m}^{\gamma, \delta}$.
The structure maps are defined by
\[
\begin{array}{rcl}
    s_{M^\mathcal{W}}(x, W_{n, m}^{\gamma, \delta})
        &=&     (x, W_n^\gamma),                    \\\\
    t_{M^\mathcal{W}}(x, W_{n, m}^{\gamma, \delta})
        &=&     (x, W_m^\delta),                    \\\\
    i_{M^\mathcal{W}}(x, W_{n, m}^{\gamma, \delta})
        &=&     (x, W_{m, n}^{\delta, \gamma}),     \\\\
    u_{M^\mathcal{W}}(x, W_n^\gamma)
        &=&     (x, W_{n, n}^{\gamma, \gamma});
\end{array}
\]
a composable pair of arrows is of the form
$\left((x, W_{t, n}^{\nu, \gamma}),
(x, W_{n, m}^{\gamma, \delta})\right)$,
and the composition is defined as
\[
    m_{M^\mathcal{W}}\left((x, W_{t, n}^{\nu, \gamma}),
    (x, W_{n, m}^{\gamma, \delta})\right)
    =
    (x, W_{t, m}^{\nu, \delta}).
\]

Define a left $\Gamma$-action on $M^\mathcal{W}$ by
\[
\begin{array}{rcl}
    \Gamma \times (M^{\mathcal W})_0
            &   \rTo &  (M^{\mathcal W} )_0        \\\\
    (\gamma^\prime, (x, W_n^\gamma))
            &   \rMapsto     &   (\gamma^\prime x, W_n^{\gamma^\prime \gamma})
\end{array}
\]
and
\[
\begin{array}{rcl}
    \Gamma \times (M^{\mathcal W} )_1
            &   \rTo &  (M^{\mathcal W} )_1        \\\\
    (\gamma^\prime, (x, W_{n, m}^{\gamma, \delta }))
            &   \rMapsto     &   (\gamma^\prime x,
            W_{n, m}^{\gamma^\prime \gamma, \gamma^\prime \delta}).
\end{array}
\]
The following proposition is straightforward.

\begin{proposition}
\label{prop-Gammaactionongrpoid}

The above is an action of the group $\Gamma$ on the Lie groupoid
$M^{\mathcal W}$.

\end{proposition}

Hence, we have the following.

\begin{definition}
\label{def-crossprodMgammaw}

Let $M^{\mathcal W}_\Gamma  = \Gamma \ltimes M^{\mathcal W}$  be
the groupoid crossed product of $\Gamma$ with the groupoid
$M^{\mathcal W}$ with respect to the above action of $\Gamma$. In
particular we have
\[
\begin{array}{rcl}
    s_{M^{\mathcal W}_\Gamma}
    \left( \gamma, (x, W_{n, m}) \right)
        &=&
    \left( x, W_n \right),
            \\\\
    t_{M^{\mathcal W}_\Gamma}
    \left( \gamma, (x, W_{n, m}) \right)
        &=&
    \left( \gamma x, W_m \right),
            \\\\
    i_{M^{\mathcal W}_\Gamma}
    \left( \gamma, ( x, W_{n, m}) \right)
        &=&
    \left( \gamma^{-1}, (\gamma x, W_{m,n})  \right),
            \\\\
    u_{M^{\mathcal W}_\Gamma}
    \left( x, W_n \right)
        &=&
    \left( 1, (x, W_{ n, n}) \right);
\end{array}
\]
a composable pair is of the form
$\left( \gamma, (x, W_{l,n}) \right), \left(\delta, (\gamma x, W_{n,m}) \right)$,
and the composition is given by
\[
    m_{M^{\mathcal W}_\Gamma}
    \left[
    \left( \gamma, (x, W_{l,n}) \right),
    \left(\delta, (\gamma x, W_{n,m}) \right)
    \right]
        =
    \left( \gamma\delta, (x, W_{l, m}) \right).
\]

\end{definition}

\begin{proposition}
\label{prop-mgammamnfld}

The groupoid  $M^{\mathcal W}_\Gamma$  is Morita equivalent to
$M_\Gamma  $ with its trivial groupoid structure.

\end{proposition}

\begin{proof}

In fact, there is a strong equivalence from $M^{\mathcal W}_\Gamma$ to
$M_\Gamma$ defined on objects by $(\gamma, (x, W_n)) \mapsto \pi_\Gamma(x)$,
and on arrows by mapping $(\gamma, (x, W_{n,m}))$ to the unit over $\pi_\Gamma(x)$.
That this map is a strong equivalence is easy to check.

\end{proof}

In the same way, one can prove the following.

\begin{proposition}
\label{prop-refinementcovers}

If ${\mathcal {\tilde W}} $ is a refinement of  ${\mathcal W}$,
then the natural groupoid morphism $\rho^{\mathcal {\tilde
W}}_{\mathcal W} : M^{\mathcal {\tilde W}}_\Gamma \to M^{\mathcal
W}_\Gamma $ is a strong equivalence.

\end{proposition}

Note that Proposition \ref{prop-mgammamnfld} implies that the Morita
equivalence class of the groupoid $M^{\mathcal{\tilde W}}_\Gamma$ does not depend on the metric used to define it.  More concretely using Proposition \ref{prop-refinementcovers}, if given two metrics on $M_\Gamma$ with corresponding covers $\mathcal{U}_1$ and $\mathcal{U}_2$ (inducing covers $\mathcal{W}_1$ and
$\mathcal{W}_2$ of $M$), one can define a strictly smaller metric
and corresponding cover $\mathcal{U}_3$ that refines both $\mathcal{U}_1$
and $\mathcal{U}_2$.

\begin{definition}
\label{def-multiloopgroupoidcover}

Let $Q$ be a smooth orbifold presented by the orbifold groupoid
$\mathcal{G}$, and let $\mathcal{W}$ be a cover of $M$ constructed
from an admissible cover of $M_\Gamma$ as above.  The {\bf
$M_\Gamma$-multiloop groupoid of $\mathcal{G}$ corresponding to
$\mathcal{W}$}  is defined to be the groupoid $\mathcal{ML(W;
G)}_{M_\Gamma}$ where
\[
    \left( \mathcal{ML(W; G)}_{M_\Gamma} \right)_0
    =
    \mbox{HOM}(M^{\mathcal { W}}_{M_\Gamma}, \mathcal{G})
\]
is the set of Lie groupoid homomorphisms from
$M_\Gamma^{\mathcal{W}}$ to $\mathcal{G}$. The arrows in
$\mathcal{ML(W;G)}_{M_\Gamma}$ are defined as follows. For any two
elements $\Phi, \Psi \in
\mbox{HOM}(M^{\mathcal{W}}_\Gamma,\mathcal{G})$, an arrow from
$\Psi$ to $\Phi$ is a map $\Lambda : (M^{\mathcal W}_\Gamma )_1 \to
(\mathcal{G})_1$ such that the following diagram commutes
\[
\begin{diagram}[width=3cm]
    (M^{\mathcal W}_\Gamma  )_1     &   \rTo^\Lambda      &       G_1       \\
    \dTo<{s_{M^{\mathcal W}_\Gamma}\times t_{M^{\mathcal W}_\Gamma}}
                                    &                    &
                                \dTo_{s\times t}\\
    (M^{\mathcal W}_\Gamma)_0 \times (M^{\mathcal W}_\Gamma)_0
                                    &   \rTo^{\Psi_0 \times \Phi_0} &        G_0 \times G_0  \\
\end{diagram}
\]
and  such that for every $(\gamma, (x, W_n)) \in \left(M^{\mathcal W}_\Gamma\right)_1$, we have
\begin{equation}
\label{eq-definelambdaloop}
\begin{array}{rcl}
    \Lambda (\gamma, (x, W_n))
    &=&
    \Psi_1 (\gamma, (x, W_n)) \, \Lambda \left[u_{M^{\mathcal W}_\Gamma}\circ s_{M^{\mathcal W}_\Gamma}
    (\gamma, (x, W_n))\right]       \\\\
    &=&
    \Lambda \left[u_{M^{\mathcal W}_\Gamma} \circ t_{M^{\mathcal W}_\Gamma} (\gamma, (x, W_n))\right]
    \, \Phi_1 (\gamma, (x, W_n))
\end{array}
\end{equation}
where as usual $\Phi_1$ and $\Psi_1$ denote the maps on arrows given by $\Phi$
and $\Psi$, respectively.  Note that
the above product is taken in $G_1$ so that the target of the right element is equal to the
source of the left.

If $\Lambda : \Psi \to \Phi$ and $\Omega :
\Phi \to \Xi$, then the composition $\Omega \circ \Lambda$ is
defined by
\[
\begin{array}{c}
    \Omega\circ \Lambda\left[
    u_{M^{\mathcal W}_\Gamma}\circ t_{M^{\mathcal W}_\Gamma}
    (\gamma, (x, W_n))\right]   \\\\
    =
    \Lambda\left[u_{M^{\mathcal W}_\Gamma}\circ t_{M^{\mathcal W}_\Gamma}
    (\gamma, (x, W_n))\right]
    \Omega \left[u_{M^{\mathcal W}_\Gamma}\circ t_{M^{\mathcal W}_\Gamma}
    (\gamma, (x, W_n))\right]
\end{array}
\]
and
\[
    \Omega\circ \Lambda (\gamma,(x, W_n))
    =
    \Omega\circ \Lambda
    \left[u_{M^{\mathcal W}_\Gamma}\circ t_{M^{\mathcal W}_\Gamma}(\gamma,(x, W_n))\right] \,
    \Xi_1 (\gamma,(x, W_n))
\]

\end{definition}

Under the compact-open topology, $\mathcal{ML(W; G)}_{M_\Gamma}$ is a
topological groupoid. Note that for each arrow in $\Lambda \in
\mathcal{ML(W;G)}_{M_\Gamma}$ from $\Phi$ to $\Psi$, $\Lambda \circ
u_{M^{\mathcal W}_{M_\Gamma}}$ is a natural transformation from $\Phi$
to $\Psi$ (see \cite[Definition 1.40]{ademleidaruan}).

Compare the following to \cite[Definition 3.2.2]{luperciouribeloop}.

\begin{definition}
\label{def-multiloopgroupoid}

Let $Q$ be a smooth orbifold  presented by the orbifold groupoid
$\mathcal{G}$.  The {\bf $M_\Gamma$-multiloop groupoid of
$\mathcal{G}$} denoted $\mathcal{ML(G})_{M_\Gamma}$ is the colimit
of the $\mathcal{ML( W;G)}_{M_\Gamma}$ over all admissible covers of $M_\Gamma$
partially ordered by inclusion of cover charts.

\end{definition}

In order to show that the $M_\Gamma$-multiloop
groupoid $\mathcal{ML(G})_{M_\Gamma}$ of an orbifold $Q$ presented by the orbifold groupoid
$\mathcal{G}$ is \'{e}tale, we first have the following.

\begin{lemma}
\label{lem-arrowsmultiloopgroupoid}

Let  $\mathcal{ML(G})_{M_\Gamma}$ be the $M_\Gamma$-multiloop
groupoid of an orbifold $Q$ presented by the orbifold groupoid $\mathcal{G}$. Then any
arrow $\Lambda: \Psi \to \Phi$ is completely determined by $\Psi$
and by $\Lambda \circ u_{M^{\mathcal W}_\Gamma}(x, W_n)$ for any
$(x, W_n) \in (M^{\mathcal W}_\Gamma)_0$.

\end{lemma}

\begin{proof}

Straightforward from the definitions; see \cite[Lemma
3.2.4]{luperciouribeloop}.   On any given $W_n$, $\Lambda$ is
determined by $\Psi$, a single value $\Lambda \circ u_{M^{\mathcal
W}_\Gamma}(x, W_n)$, and Equation \ref{eq-definelambdaloop}.  This
determines $\Lambda$ on each chart $W_m$ such that $W_n \cap W_m \neq \emptyset$,
and hence recursively on on every chart.

\end{proof}

\begin{proposition}
\label{prop-multiloopetale}

Let  $\mathcal{ML(G})_{M_\Gamma}$ be the $M_\Gamma$-multiloop
groupoid of an orbifold $Q$ presented by the orbifold groupoid $\mathcal{G}$. Then
$\mathcal{ML(G})_{M_\Gamma}$ is \'{e}tale.

\end{proposition}

\begin{proof}

Again straightforward from the definitions.  An arrow $\Lambda: \Psi
\to \Phi$ is determined by $\Psi$ and a single value $\Lambda \circ
u_{M^{\mathcal W}_\Gamma}(x, W_n)$. The result then follows from the
fact that the isotropy groups of  $\mathcal{G}$ are finite.

\end{proof}

The next result follows \cite[Section 3.4]{luperciouribeloop}.

\begin{proposition}
\label{prop-grpoidconstants}

Let  $\mathcal{ML(G}_i)_{M_\Gamma}$ be the $M_\Gamma$-multiloop
groupoid of orbifolds $Q_i$ presented by the orbifold groupoids $\mathcal{G}_i$ for
$i=1,2$. A groupoid homomorphism $e: Q_1 \rightarrow Q_2$ induces a
homomorphism of $M_\Gamma$-multiloop groupoids $e_\mathcal{ML}:
\mathcal{ML(G}_1)_{M_\Gamma} \rightarrow \mathcal{ML(G}_2)_{M_\Gamma}$.  If
$e$ is a strong equivalence, then $e_\mathcal{ML}$ is a strong
equivalence.

\end{proposition}

\begin{proof}

The proof is identical to that of
\cite[Definition 3.4.1 and Lemma 3.4.2]{luperciouribeloop}.

\end{proof}


\subsection{The $M_\Gamma$-Multiloops when $\Gamma$ is a
Subgroup of a Contractible Abelian Group}
\label{subsec-gammasubgroupcontractabel}

In this subsection, we assume that $\Gamma$ is a subgroup of a
contractible abelian Lie group $T$; in the notation of Subsection
\ref{subsec-intromultiloop}, $T=M$ and $T/\Gamma = M_\Gamma$.
Following \cite[Section 3.6]{luperciouribeloop}, we recover the
$\Gamma$-sectors of $Q$ from the fixed points of the $T/\Gamma$-multiloop
groupoid. First, we define a $T$-action on the
$T/\Gamma$-multiloop groupoid $\mathcal{ML(G})_{T/\Gamma}$.

\begin{definition}
\label{def-Taction}

Suppose $\Gamma$ is a subgroup of a contractible abelian Lie group $T$.
Let $\mathcal{ML(W;G)}_{T/\Gamma}$ be the $T/\Gamma$-multiloop
groupoid associated to the cover $\mathcal{W}$ of $T$ given by
Definition \ref{def-multiloopgroupoidcover}.
For each $W_n = \coprod_{\gamma \in \Gamma} W_n^\gamma$ and $t \in
T$, let $W_n^t$ denote the translate $tW_n = \{ tx : x \in W_n \}$,
and let $\mathcal{W}^t$ denote the translated cover $\{ W_n^t
\}_{n\in\mathbb{N}}$.  Note that this introduces no ambiguity;
$W_n^t$ has the same meaning as in Subsection
\ref{subsec-intromultiloop} when $t \in \Gamma \leq T$. Then $T$ acts on
$\coprod_{t \in T} \mathcal{W}^t$ via $(s, (x, W_n^t)) \mapsto
(sx,W_n^{st})$ for $s \in T$.  As $T$ is abelian, this action
descends to a $T$-action on the cover $\coprod_{n \in \mathbb{N},
t\Gamma \in T/\Gamma} U_n^{t\Gamma}$ of $T/\Gamma$-translates of
$\mathcal{U}$ in the same way.

Now define an action of $T$ on $\coprod_{t \in T}
\mathcal{ML(W}^t;\mathcal{G)}_{T/\Gamma}$ by
\[
\begin{array}{rcl}
    T \times \left(\coprod\limits_{t \in T} \mathcal{ML(W}^t;\mathcal{G)}_{T/\Gamma}\right)_0
    &\rTo&
    \left(\coprod\limits_{t \in T} \mathcal{ML(W}^t;\mathcal{G)}_{T/\Gamma}\right)_0       \\\\
    (t, \Psi)
    &\rMapsto&           \Psi^t
\end{array}
\]
where $\Psi^t$ is defined by
\[
\begin{array}{rcl}
    \Psi_0^t(x, W_n^t )
    &=&
    \Psi (t^{-1}x, W_n),     \\\\
    \Psi_1^t(\gamma, (x, W_{n,m}^t))
    &=&
    \Psi_1 (\gamma, (t^{-1}x, W_{n,m})).
\end{array}
\]
Taking the colimit, we obtain an action of $T$ on
$\mathcal{ML(G})_{T/\Gamma}$.

\end{definition}

Now consider the subgroupoid $\mathcal{ML(G})_{T/\Gamma}^T$ of
$\mathcal{ML(G})_{T/\Gamma}$ consisting of elements fixed by the
action of $T$.  In Theorem \ref{thrm-loopmoritatosectors} below, we
will show that $\mathcal{ML(G})_{T/\Gamma}^T$ is Morita equivalent to
$\mathcal{G}^\Gamma$, the groupoid presenting the $\Gamma$-sectors
of $Q$.  When $\Gamma = \Z$, this coincides with \cite[Theorem 3.6.4
and Proposition 3.6.6]{luperciouribeloop}; see also
\cite{defernetluperciouribe}. The following two lemmas can be proved
in the same way as in \cite[Lemmas 3.6.2 and
3.6.3]{luperciouribeloop}. We give the proof of Lemma
\ref{lem-lifttotrivcover} explicitly, because it is important for
the proof of Theorem \ref{thrm-loopmoritatosectors}.

\begin{lemma}
\label{lem-Psilocconstant}

For any object $\Psi$ of $\mathcal{ML(G})_{T/\Gamma}^T$, $\Psi_0$ and
$\Psi_1$ are locally constant.

\end{lemma}

\begin{lemma}
\label{lem-lifttotrivcover}

For any object  $\Psi$ of $\mathcal{ML(G})_{T/\Gamma}^T$, there is
another object $\Phi$ of $\mathcal{ML(G})_{T/\Gamma}^T$ defined over the
trivial cover of $T$ by one chart such that there is an arrow
$\Lambda $ connecting $\Psi$ and $\Phi$.

\end{lemma}

\begin{proof}

Given a point $(x, W) \in (T^\mathcal{W})_0$, if $(y,W^\prime)\in
(T^\mathcal{W})_0$ is another point in $(T^\mathcal{W})_0$, then
there exists a finite collection $W^1, \dots, W^r$ of sets in the
cover $\mathcal{W}$ such that $W^1 = W$, $W^r = W^\prime$, and $W^i
\cap W^{i+1} \neq \emptyset$ in $T$ for $i = 1, 2, \ldots , r - 1$.

For each $(z_i, W^i)$ such that $z_i \in W^i \cap W^{i+1}$, we have that
$\Psi_1(1, (z_i, W^i \cap W^{i+1}))$ is an arrow from $\Psi_0(z_i,
W^i)$ to $\Psi_0(z_i, W^{i+1})$.  We define the arrow
$A_i^{i+1}(z_i) \in (T^\mathcal{W})_1$ from $\Psi_0(z_i,W^i)$ to
$\Psi_0(z_i,W^{i+1})$ by
\[
    A_i^{i+1}(z_i) = \Psi_1(1, (z_i, W^i \cap W^{i+1}))
\]
and note that as $\Psi_0$ and $\Psi_1$ are locally constant by Lemma
\ref{lem-Psilocconstant}, $A_i^{i+1}$ does not depend on $z_i$.
Choosing one $z_i \in W^i \cap W^{i+1}$ for each $i$, we set
\[
    A_W^{W^\prime}
    =
    \Psi_1(1, (z_{r-1}, W^{r-1} \cap W^r)) \,
    \Psi_1(1, (z_{r-2}, W^{r-2} \cap W^{r - 1})) \,
    \cdots \,
    \Psi_1(1, (z_1, W^1 \cap W^2)).
\]
Note that the definition of $A_W^{W^\prime}$ depends only on the
sets $W$ and $W^\prime$, and that $(A_W^{W^\prime})^{-1} =
A_{W^\prime}^W$.

Define a morphism $\Phi : T^\mathcal{W}_\Gamma \to \mathcal{G}$ by
\[
\begin{array}{rcl}
    \Phi_0(z, W_n)   &=&     \Psi_0(x, W)           \\\\
    \Phi_1(\gamma, (z, W_{n,m}))
                &=&
        A_{W_n}^{W_m} \Psi_1(\gamma, (z, W_{m,n}))
        A_{W_n}^{W_m}.
\end{array}
\]
Because $\Psi_1$ is locally constant,
\[
    \Phi_1 (1, (z, W_{n, m}))
    =
    \Psi_1 (1, (z, W_{n, m})),
\]
and
\[
    \Phi_1 (\gamma, (z, W_{n, m}))
    =
    \Psi_1 (\gamma, (z, W_{n, m}))
\]
for each $z \in W_{n,m}$.  Hence we can define $\Phi$ on the trivial
cover of $T$ consisting of points $(z, T)$ by
\[
\begin{array}{rcl}
    \Phi_0 (z, T)   &=&     \Phi_0 (x, W)
            \\\\
    \Phi_1 (\gamma, (z, T)) &=&  \Phi_1 (\gamma, (z, W_{n,m}))
\end{array}
\]
whenever $z \in W_n \cap W_m$. Lastly, we  define an arrow $\Lambda
: \Psi \to \Phi$ by $\Lambda (1, (z, W_{n, m}) )= A_{W_n}^{W_m}$.

\end{proof}

We now have the following.

\begin{theorem}
\label{thrm-loopmoritatosectors}

There is a strong equivalence from $\mathcal{ML(G})_{T/\Gamma}^T$ to the
groupoid $\mathcal{G}^\Gamma$ of $\Gamma$-sectors of $Q$.

\end{theorem}

\begin{proof}

See the proof of \cite[Theorem 3.6.4]{luperciouribeloop}.
Given $\Psi \in \mathcal{ML(G})_{T/\Gamma}^T$, let $\Phi$ be as in the
proof of Lemma \ref{lem-lifttotrivcover}.  As $\Phi_0$ is locally
constant, $\Phi_0(y,T) = \Phi_0(1,T)$ for each $y \in T$. We have
that
\[
\begin{array}{rcl}
    s\circ \Phi_1(\gamma, (y, T))
        &=&     \Phi_0(y,T)             \\\\
        &=&     \Phi_0(\gamma y, T)     \\\\
        &=&     t\circ \Phi_1(\gamma, (y, T))
\end{array}
\]
so that each $\Phi_1(\gamma, (y, T))$ is an element of the isotropy
group $G_{\Phi_0(y)} = G_{\Phi_0(1)}$.  Hence, we can define a
homomorphism $\phi : \Gamma \to G_{\Phi_0(1)}$ by
\[
    \phi(\gamma) =  \Phi_1(\gamma, (1, T)).
\]
Clearly, the correspondence $\Psi \mapsto \phi$ is surjective, as
given any $\phi_x : \Gamma \rightarrow G_x$, one can define a $\Psi
\in \mathcal{ML({G}})_{T/\Gamma}^T$ with $\Psi_0(y,T) = x$ and
$\Psi_1(\gamma, (y,T)) = \phi_x(\gamma)$.  That this correspondence
is a strong equivalence of groupoids is straightforward.

\end{proof}

It follows that $\mathcal{ML(G})_{T/\Gamma}^T$ is Morita equivalent to
$\mathcal{G}^\Gamma$.

Note that in \cite[Proposition 3.5.3]{chenorbispace}, Chen proves that
the in the case of a proper \'{e}tale topological groupoid $\mathcal{G}$ representing
an orbispace such that the space of objects is a $T_1$ space, there is an identification
similar to that given by Theorem \ref{thrm-loopmoritatosectors} on the level of orbispaces for the
case $\Gamma = \Z$.  Using exactly the same proof with the definitions given above and
Equation \ref{eq-orbispacesectors}, we have the following.

\begin{proposition}
\label{prop-chenorbispacegen}

Let $X$ be an \'{e}tale proper orbispace, that is, an orbispace
represented by the \'{e}tale proper groupoid $\mathcal{G}$, such that
$G_0$ is a $T_1$ space.  Let $\Gamma$ be a discrete subgroup of a contractible abelian Lie group $T$.
Then the orbit space $|\mathcal{ML(G})_{T/\Gamma}^T|$ of the groupoid
$\mathcal{ML(G})_{T/\Gamma}^T$ is homeomorphic to $|\mathcal{G}^\Gamma|$.

\end{proposition}


\subsection{The $M_\Gamma$-Multiloop and
the $\Gamma$-Sectors in the General Case.}
\label{subsec-generalcase}

In the general case of $M_\Gamma$ an arbitrary manifold with
fundamental group $\Gamma$ and universal cover $\Gamma$, we have a
correspondence similar to Theorem \ref{thrm-loopmoritatosectors}. In
this case, we use the groupoid of constants, a subgroupoid of
$\mathcal{ML(G})_{M_\Gamma}$.

\begin{definition}
\label{def-grpoidconstantsgeneral}

Let $Q$ be an orbifold presented by the orbifold groupoid
$\mathcal{G}$.  The {\bf groupoid of constants}
$\mathcal{C(G})_{M_\Gamma}$ of $\mathcal{ML(G})_{M_\Gamma}$ is defined to be
the subgroupoid of $\mathcal{ML(G})_{M_\Gamma}$ consisting of the $\Phi$
such that $\sigma \circ \Phi$ is constant.  Recall that $\sigma :
\mathcal{G} \rightarrow |\mathcal{G}|$ denotes the quotient map
onto the orbit space of $\mathcal{G}$.

\end{definition}

\begin{theorem}
\label{thrm-loopmoritatosectorsgeneral}

There is a strong equivalence from the groupoid of constants
$\mathcal{C(G})_{M_\Gamma}$ to the groupoid $\mathcal{G}^\Gamma$ of
$\Gamma$-sectors of $Q$.

\end{theorem}

The proof is identical to that of Theorem
\ref{thrm-loopmoritatosectors}.

\subsection{The $M_\Gamma$-Multiloops of a Quotient
Orbifold}
\label{subsec-loopquot}

In this subsection, we specialize to the case where $M_\Gamma$ is compact and
$Q$ is presented by as the quotient of a smooth connected manifold $X$ by a compact Lie group $G$
acting locally freely (i.e. properly with discrete stabilizers).
In the case of $G$ finite, a very explicit characterization of the loop space is given in
\cite[Section 4.1]{luperciouribeloop}.  It is shown that it is
enough to consider only the homomorphims defined on the trivial cover.
Here, we briefly explain how this characterization extends readily to the case of $G$ compact
and the $M_\Gamma$-multiloops.  Throughout this section, we let $\mathcal{G}$ denote an orbifold groupoid
Morita equivalent to $G \ltimes X$.  In particular, we can take $\mathcal{G}$ to be given by
a collection of slices for the $G$-action as in Theorem \ref{thrm-localglobaldiffeogenerallie}.

First, we note the following.
The proof is similar to that of Lemma \ref{lem-lifttotrivcover} and,
given the modifications outlined in Definitions \ref{def-crossprodMgammaw},
\ref{def-multiloopgroupoidcover}, and
\ref{def-grpoidconstantsgeneral} and the local structure of quotient orbifolds
demonstrated by Theorem \ref{thrm-localglobaldiffeogenerallie},
identical to that of \cite[Lemma 4.1.1]{luperciouribeloop}.

\begin{proposition}
\label{prop-lifttotrivcovergeneral}

Let $Q$ be a quotient orbifold presented by $G \ltimes X$ with $G$ a compact Lie group acting
locally freely on the smooth manifold $X$.
Let  $\mathcal{ML(G})_{M_\Gamma}$ be the
$M_\Gamma$-multiloops for $M_\Gamma$ compact.  Then for  any
morphism $\Psi: M_\Gamma^\mathcal{W} \to \mathcal{G}$, there is a
morphism $\Phi: M_\Gamma^{ \{ M \}} \to \mathcal{G}$ subordinate to
the trivial cover of $M_\Gamma$ by $M$ and an arrow connecting
$\Psi$ to $\Phi$.

\end{proposition}

Similarly, the morphisms $\Psi$ subordinated to
the trivial cover of $M_\Gamma$ is determined by the image of
$\coprod_{i=1}^s M\times \{ \gamma_i \}$ under  $\Psi_1$
where $\{ \gamma_1, \dots , \gamma_s \}$ is  a set of generators of $\Gamma$; compare
\cite[Section 3.3]{luperciouribeloop}.

\begin{lemma}
\label{lem-morphismsdeterminedbygensgeneral}

Let $Q$ be a quotient orbifold presented by $G \ltimes X$ with $G$ a compact Lie group acting
locally freely on the smooth manifold $X$.
Every morphism in $\Psi \in \mathcal{ML(G})_{M_\Gamma}$
subordinated to the trivial cover of $M_\Gamma$ is determined
by the image of $\coprod_{i=1}^n M\times \{ g_i \}$ under  $\Psi_1$,
where $\{ g_{1}, \dots , g_{n} \}$ is  a set of generators of $\Gamma$.

\end{lemma}

\begin{proof}

Pick a set of generators $\{ \gamma_1, \ldots , \gamma_s \}$ of $\Gamma$ and
let $(x, \gamma) \in \Gamma \ltimes M$.  Then if
$\gamma = \gamma_{\alpha_1}^{\beta_1}\dots \gamma_{\alpha_s}^{\beta_s}$
is an expression of $\gamma$ in terms of these generators,
\[
    (x,\gamma)
    =
    (x, g_{\alpha_1})
    (g_{\alpha_1}x, g_{\alpha_1}^{\beta_1-1 })
    \cdots (g_{\alpha_s}^{\beta_s-1} \cdots g_{\alpha_1}^{\beta_1} x, g_{\alpha_s}).
\]
It follows that
\[
    \Psi_1(x,\gamma)
    =
    \Psi_1(x, \gamma_{\alpha_1}) \Psi_1(x\, \gamma_{\alpha_1},
    \gamma_{\alpha_1}^{\beta_1-1 }) \cdots \Psi_1 (x\, \gamma_{\alpha_1}^{\beta_1} \cdots \gamma_{\alpha_n}^{\beta_n-1}, \gamma_{\alpha_n}).
\]

\end{proof}

Fixing a generating set
$\{ \gamma_1, \ldots , \gamma_s \}$ of $\Gamma$,
it follows that there is a bijective
correspondence between the morphisms $\Psi$
subordinated to the trivial cover of $M_\Gamma$  in
$\mathcal{ML}(\mathcal{G})_\Gamma$
and the collection of pairs $(f , {\Theta})$
where $\Theta = \{ g_1, \dots , g_n \}$ is an $s$-tuple of
elements of $G$ satisfying the same relations as the $\gamma_i$,
and $f : M \rightarrow X$ is a smooth map such that
$g_i f(x) = f(\gamma_i x) $ for each $i = 1, \ldots s$.  Let $\mathcal{P}_\Theta$
denote the set of all such pairs.  Similarly,
let $\Lambda$ be an arrow between homomorphisms $\Psi = (f, {\Theta})$
and $\Phi =(f^\prime, {\Upsilon})$ with
${\Theta}= \{ g_1, \ldots , g_s \}$ and ${\Upsilon}= \{ k_1, \ldots , k_s \}$.
Using  the fact that $X$ is connected,
there is an $h \in G$ such that $k_i = h g_i h^{-1}$ and $hf (x)= f^\prime (x)$
for each $i = 1, \ldots s$ and $x \in M$.  Thus we have the following consequence
of Lemma \ref{lem-arrowsmultiloopgroupoid}.

\begin{proposition}
\label{prop-crossprodequivmultiloops}

Let $G$ act on ${\mathcal P}_\Theta$ via
\[
    [h, (f , {\Theta})]  \rMapsto (hf , h \Theta h^{-1}),
\]
where $h {\Theta}h^{-1}$ indicates pointwise conjugation of the $s$-tuple
$\Theta$.  The crossed product groupoid $G \ltimes {\mathcal P}_\Theta$ is
Morita equivalent to the orbifold $M_\Gamma$-multiloop groupoid
$\mathcal{ML(G})_{M^\Gamma}$ of $Q$.

\end{proposition}

For each $s$-tuple $\Theta = (g_1, \ldots, g_s)$, let $C_G(\Theta)$ denote the
centralizer of the subgroup generated by the $g_i$.  Techniques identical to those in
Lemma \ref{lem-globalsectorglobalactiondef} demonstrate that the crossed-product
groupoid $G \ltimes {\mathcal P}_\Theta$ is given Morita equivalent to
\[
    \coprod\limits_{(\Theta)} \left(  C_G (\Theta ) \ltimes {\mathcal P}_\Theta  \right),
\]
where the sum is over the $G$-conjugacy classes $(\Theta)$ of the $s$-tuples $\Theta$.  Hence we have
the following.

\begin{corollary}
\label{cor-loopstructureequivmultiloops}

The groupoid
\[
    \coprod\limits_{(\Theta)} \left(  C_G (\Theta ) \ltimes {\mathcal P}_\Theta  \right),
\]
is Morita equivalent to the orbifold $M_\Gamma$-multiloop groupoid
$\mathcal{ML(G})_{M^\Gamma}$ of $Q$.

\end{corollary}

Note that in \cite[Equation 2-12, page 808]{tamanoi2}, the multiloop $\Gamma$-sectors for
global quotients are given by
\[
    \mathbb{L}_{M_\Gamma} (X;G)
    =
    \coprod_{(\phi)} \mbox{Map}_{\phi}(M; X) /C_G(\phi)
\]
where the space $\mbox{Map}_{\phi}$ is defined as
\[
    \mbox{Map}_{\phi}(M; X)
    =
    \{ f : M \to X | f(\gamma x)
    = \phi(\gamma) f(x) \;\forall x \in M, \gamma \in \Gamma  \}
\]
(with notation modified to our case, and the $\Gamma$-action on $M$ expressed as a left action).
Hence $\mathbb{L}_{M_\Gamma}$ coincides with the groupoid in Corollary \ref{cor-loopstructureequivmultiloops} in the case of $G$ finite.

Restricting to the groupoid of constant maps $(f, \Theta)$ in ${\mathcal ML (G)}_{M_\Gamma}$
as in Definition \ref{def-grpoidconstantsgeneral}, we have that if $\Theta = (g_1, \ldots, g_s)$, then
\[
     g_i f(x) = f(\gamma_i x) = f(x)
\]
for each $x \in X$ and $i = 1, \ldots , s$.  Hence the image of $f$ is fixed by each $g_i$.

\begin{corollary}

The subgroupoid of constants ${\mathcal CL (\mathcal{G})}_{M_\Gamma}$ of
${\mathcal ML (G)}_\Gamma $ is given by
\[
    \coprod\limits_{(\Theta )}
    \left(  C_G (\Theta ) \ltimes X^{\langle \Theta \rangle}\right),
\]
where the sum is over the $G$-conjugacy classes $(\Theta)$ of the $s$-tuples $\Theta$,
and is Morita equivalent to the groupoid $\mathcal{G}^\Gamma$ presenting the $\Gamma$-sectors $\tilde{Q}_\Gamma$.

\end{corollary}


\bibliographystyle{amsplain}

\end{document}